\documentclass[12pt]{amsart}

\usepackage{epsfig}
\usepackage{amsmath}
\usepackage{amssymb}
\usepackage{graphicx}

\usepackage{tikz}

\usepackage{hyperref}
\usepackage{a4wide}

\newtheorem{propo}{Proposition}[section]
\newtheorem{defi}[propo]{Definition}

\newtheorem{lemma}[propo]{Lemma}
\newtheorem{corol}[propo]{Corollary}

\newtheorem{theo}[propo]{Theorem}

\newtheorem{prob}[propo]{Problem}

\newtheorem{rem}[propo]{Remark}

\newcommand{\bl}{\begin{lemma}}
\newcommand{\el}{\end{lemma}}

\newcommand{\Inc}{\mathop{\rm Inc}\nolimits}

\def\d12{{_{12}}}

\usepackage{indentfirst,latexsym,bm}

\def\PSL{{\rm PSL}}

\def\Aut{{\rm Aut}}

\def\diam{{\rm diam}}

\def\K{{\rm K}}

\def\H{{\rm H}}

\begin{document}
\title{Finite 2-arc-transitive strongly regular graphs and 3-geodesic-transitive  graphs
}

\thanks{The first author acknowledges the hospitality of the Centre
for the Mathematics of Symmetry and Computation of UWA, where this research was carried
out.}
\thanks{The first author is supported by NSFC (11661039) and NSF of Jiangxi (2018ACB21001, 20171BAB201010,20171BCB23046,GJJ170321) and the second author acknowledges the Australian Research Council grant DP130100106.}

\author[W. Jin]{Wei Jin}
\address{Wei Jin\\School of Statistics\\
Jiangxi University of Finance and Economics\\
 Nanchang, Jiangxi, 330013, P.R.China}
 \address{Research Center of Applied Statistics\\
Jiangxi University of Finance and Economics\\
 Nanchang, Jiangxi, 330013, P.R.China}
\email{jinweipei82@163.com}
\author[C. E. Praeger]{Cheryl E. Praeger}
\address{Cheryl E. Praeger\\CMSC and Department  of Mathematics and Statistics\\
The University of Western Australia\\
Crawley, WA 6009, Australia} \email{cheryl.praeger@uwa.edu.au}


\maketitle

\begin{abstract}
We classify all the $2$-arc-transitive strongly regular graphs,  and use this classification
to study the family of finite $(G,3)$-geodesic-transitive graphs of girth $4$ or $5$  for some group $G$ of automorphisms.
For this application we first give a reduction result on the latter family of graphs: let $N$ be a  normal subgroup of
$G$ which has at least $3$ orbits on vertices. We show that
$\Gamma$ is a cover of its  quotient  $\Gamma_N$ modulo the $N$-orbits, and that either
$\Gamma_N$ is  $(G/N,3)$-geodesic-transitive of the same girth as $\Gamma$, or
$\Gamma_N$ is a $(G/N,2)$-arc-transitive strongly regular graph, or
$\Gamma_N$ is a complete  graph with $G/N$ acting 3-transitively on vertices.
The  classification of   $2$-arc-transitive strongly regular graphs allows us to  characterise the $(G,3)$-geodesic-transitive covers $\Gamma$
when $\Gamma_N$ is complete or strongly regular.

\end{abstract}

\vspace{2mm}

 \hspace{-17pt}{\bf Keywords:}  $3$-geodesic-transitive graph, strongly regular graph, automorphism group.

 \hspace{-17pt}{\bf Math. Subj. Class.:} 05E18; 20B25

\section{Introduction}

A \emph{geodesic} from a vertex $u$ to a vertex $v$ in a graph $\Gamma$ is
a path  of the shortest length  from $u$ to $v$ in $\Gamma$, and  is
called an \emph{$s$-geodesic} if the distance between $u$ and $v$ is $s$.
Then $\Gamma$ is said to be \emph{$(G,s)$-geodesic-transitive} if it has an $s$-geodesic, and for each   $i\leq s$, the  automorphism group
$G$ is transitive on the set of  $i$-geodesics of $\Gamma$. The systematic investigation of $s$-geodesic-transitive graphs was initiated recently. The  possible local structures of  2-geodesic-transitive graphs were determined in \cite{DJLP-clique}. Then
Devillers, Li  and the authors
\cite{DJLP-2} classified $2$-geodesic-transitive  graphs of valency 4. Later, in \cite{DJLP-cayleyred1},  a
reduction theorem for the family of normal
$2$-geodesic-transitive Cayley graphs  was proved and
those which are complete multipartite graphs were also classified.
Our focus in this paper is on 3-geodesic-transitive graphs.

For a positive integer $s$, an \emph{$s$-arc} of $\Gamma$ is a sequence of
vertices $(v_0,v_1,\ldots,v_s)$  in  $\Gamma$ such that
$v_i,v_{i+1}$ are adjacent and $v_{j-1}\neq v_{j+1}$ where $0\leq
i\leq s-1$ and $1\leq j\leq s-1$. In particular, 1-arcs are called
\emph{arcs}.  The graph $\Gamma$ is said to be \emph{$(G,s)$-arc-transitive}
if, for each $i\leq s$, the automorphism group
$G$ is transitive on the set of  $i$-arcs of $\Gamma$.  Study of  $s$-arc-transitive graphs originates
from Tutte  \cite{Tutte-1,Tutte-2},
who proved that there are no $6$-arc-transitive
cubic graphs, and for such graphs the order of the stabiliser  of a vertex  is at
most 48. This seminal result stimulated greatly the study of $s$-arc-transitive graphs.
About twenty years later, relying on the
classification of finite 2-transitive groups (which in turn depends on the  finite
simple group classification), Weiss \cite{weiss} proved that there
were no $8$-arc-transitive graphs with valency at least three.
Moreover, for each $s\leq 5$ and $s=7$, $s$-arc-transitive graphs exist which are not $(s+1)$-arc-transitive, but there are no such graphs for $s=6$.
Many other results have been proved for $s$-arc-transitive graphs,
 see
\cite{ACMM-1996,IP-1,Li-abeliancay-2008}.

On the other hand, there is no upper bound on $s$ for $s$-geodesic-transitivity \cite[Theorem 1.1]{WJ-2015au}.
Clearly, every $s$-geodesic is an $s$-arc, but some $s$-arcs
may not be $s$-geodesics, even for small values of $s$.
If $\Gamma$ has girth 3 (the \emph{girth} of $\Gamma$ is
the length of the shortest cycle in $\Gamma$), then  2-arcs contained in $3$-cycles are
not 2-geodesics.
If $\Gamma$ has girth 4 or 5, then 3-arcs contained in $4$-cycles or $5$-cycles are
not 3-geodesics. The graph in Figure \ref{fig1} is the Hamming  graph $\H(3,2)$, and is
$(G,3)$-geodesic-transitive but not $(G,3)$-arc-transitive with
valency 3 and girth 4 where $G$ is the full automorphism group. Thus  the family of  $(G,3)$-arc-transitive
graphs is properly contained in the family of $(G,3)$-geodesic-transitive graphs.
We study $s$-geodesic-transitive graphs that  are  not $s$-arc-transitive for $s=3$. This problem was studied earlier for the case $s=2$,
refer to \cite{DJLP-clique,DJLP-2,DJLP-cayleyred1,DJLP-prime}.
For $s=3$, the valency 4 examples have been classified in
\cite{WJ-2015au}, where it is shown also that there are  examples with
unboundedly  large diameter and valency.
In this paper, we introduce a general framework
for describing all the graphs with these properties.

\begin{figure}[t]\label{fig1}
\centering

\begin{tikzpicture}

\draw (1,1)-- (3,0); \draw (1,1)-- (3,1);

\draw (1,1)-- (3,2);

\draw (3,2)-- (5,2);

\draw (3,1)-- (5,2); \draw (3,1)-- (5,0);

\draw (3,0)-- (5,1); \draw (3,0)-- (5,0);



\draw (3,2)-- (5,1);

\draw (7,1)-- (5,0); \draw (7,1)-- (5,1);

\draw (7,1)-- (5,2);


\filldraw[black] (1,1) circle (2pt)  (5,0) circle (2pt);

\filldraw[black] (3,0) circle (2pt) (3,1) circle (2pt);

\filldraw[black] (3,2) circle (2pt);

\filldraw[black] (5,2) circle (2pt)  (5,1) circle (2pt);

\filldraw[black] (7,1) circle (2pt);

\end{tikzpicture}
\caption{$\H(3,2)$}
\end{figure}

We study normal quotients. Let $\Gamma$ be a $G$-vertex-transitive
graph. If $N$ is a vertex-intransitive normal subgroup of $G$, then  the \emph{quotient graph}
$\Gamma_{N}$ of $\Gamma$  is  the graph whose vertex set is the set of $N$-orbits,
such that two $N$-orbits  $B_i,B_j$ are adjacent in $\Gamma_{N}$ if
and only if there exist $x\in B_i, y\in B_j$ such that $x,y$ are
adjacent in $\Gamma$. Such quotients $\Gamma_N$ are often referred to as
\emph{$G$-normal quotients} of $\Gamma$ relative to $N$.
Sometimes    $\Gamma$ is  a
\emph{cover} of $\Gamma_{N}$, that is to say, for each edge $\{B_i,B_j\}$ of
$\Gamma_{N}$ and $v\in B_i$, $v$ is adjacent to exactly one vertex in $B_j$.
In this case we say that $\Gamma$ is  a
\emph{$G$-normal cover} of $\Gamma_N$ relative to $N$.

A   connected regular graph  is said to be \emph{strongly regular} with parameters
$(n,k,a,c)$ if it has valency $k$, vertex set of size $n$, every pair of adjacent vertices
has $a$ common neighbours, and every pair of distinct non-adjacent vertices has $c$ common neighbours.

\medskip
Our first  theorem is   a  reduction result on the family of  $3$-geodesic-transitive graphs of  girth 4 or 5.
It  describes the various possibilities for the girth and diameter of normal quotients.
Note that $3$-geodesic-transitive  graphs have diameter at least 3 and are therefore not
$3$-arc-transitive.

\begin{theo}\label{3gt-redtheo-1}
Let $\Gamma$ be a connected  $(G,3)$-geodesic-transitive  graph of girth $4$ or $5$.
Let $N$ be a  normal subgroup of $G$ with at least $3$ orbits on the vertex set. Then $\Gamma$ is a cover of $\Gamma_N$,   $\Gamma_N$ is $(G/N,s')$-geodesic-transitive where
$s'=\min\{3,\diam(\Gamma_N)\}$, and one of the following holds:

\begin{itemize}
\item[(1)] $\Gamma_N$ is a complete graph.

\item[(2)] $\Gamma_N$ is a $(G/N,2)$-arc-transitive strongly regular graph with  girth $4$ or $5$.


\item[(3)] $\Gamma_N$  has diameter at least $3$ and the same girth as $\Gamma$.

\end{itemize}

\end{theo}

The normal quotient graphs in Theorem \ref{3gt-redtheo-1} (2) are
$2$-arc-transitive strongly regular  graphs. In order to  classify $3$-geodesic-transitive  graphs of girth $4$ or $5$,
we need to know the $2$-arc-transitive  strongly regular graphs, and our second theorem determines all such  graphs.

\begin{theo}\label{diam2-theo-1}
Let $\Gamma$ be a $2$-arc-transitive  strongly regular graph.  Then either

\begin{itemize}
\item[(1)] $\Gamma$ has girth $4$ and is one of the following graphs: $\K_{m,m}$ with $m\geq 2$, the Higman-Sims graph,  the Gewirtz graph, the $M_{22}$-graph, or the folded $5$-cube $\Box_5$; or

\item[(2)] $\Gamma$ has girth $5$ and is one of the following graphs: $C_5$,
the Petersen graph,  or the Hoffman-Singleton graph.
\end{itemize}

\end{theo}

\medskip

\begin{table}[h]\caption{Quotient graphs $\Gamma_N$ for Theorem \ref{girth45-redtheo-1}}\label{table-cover-diam2}
\begin{tabular}{lllll}
\hline
 $\Gamma_N$ &  $\Gamma$  & Condition & Ref \\
\hline\hline
$\K_r$ &  $\K_{r,r}-r\K_2$, or   &   $r\geq 3$ & Lemma \ref{quot-comp-2} \\
  & $[HoS]_2 $       & &  \\
\hline
$\Box_5$ & $\H(5,2)$, or  &   & Lemma \ref{cover-petersen} \\
& Armanios-Wells graph    & &  \\
\hline
Petersen graph    & Dodecahedron &   & Lemma \ref{cover-petersen} \\
\hline
HiS    & SDC of $\Gamma_N$  &    & Lemma \ref{cover-hs-1} \\
 \hline
Gewirtz graph  &  SDC  of $\Gamma_N$   &      & Lemma \ref{cover-g-1} \\
 \hline
$M_{22}$ graph  &  SDC  of $\Gamma_N$   &      & Lemma \ref{cover-m22-c3} \\
\hline
$\K_{r,r}$   &   Hadamard graph, or    &   $m=2$, $r\geq 3$ & Lemma \ref{cover-bipart-1} \\
& $m\K_{r,r}$ from $RGD(r,c_2,m)$          &  $m>2$, $r=m.c_2$  &  \\
\hline
\end{tabular}

\end{table}

\begin{table}[h]\caption{Non-$(G,4)$-distance-transitive graphs $\Gamma$ for Theorem \ref{girth45-redtheo-1}}\label{table-cover-diam}
\begin{tabular}{llll}
\hline
 $\Gamma_N$ & $G_u$ &    Ref \\
\hline\hline
 Gewirtz graph  &   $PGL(2,9)$ or $P\Gamma L(2,9)$       & Lemma \ref{cover-g-1} \\
 \hline
 $M_{22}$ graph  &    $\mathbb{Z}_2^4:A_6$ or $\mathbb{Z}_2^4:S_6$       & Lemma \ref{cover-m22-c3} \\
\hline
$\K_{r,r}$, $r\geq 3$   &     &     Lemma \ref{cover-bipart-1} \\
\hline
\end{tabular}

\end{table}


Our third theorem,  using the  classification of   $2$-arc-transitive strongly regular graphs in Theorem \ref{diam2-theo-1},   characterises all the $(G,3)$-geodesic-transitive covers $\Gamma$
when $\Gamma_N$ is in part (1) or (2) of Theorem \ref{3gt-redtheo-1}.
For each possible quotient $\Gamma_N$, we obtain all normal covers $\Gamma$ explicitly, except for three cases for $\Gamma_N$ where we can only classify the graphs  when $\Gamma$ is $(G,4)$-distance-transitive (see Table \ref{table-cover-diam}).

\begin{theo}\label{girth45-redtheo-1}
Let $\Gamma$ be a connected  $(G,3)$-geodesic-transitive  graph of girth  $4$ or $5$.
Suppose that $G$ has a normal subgroup $N$ such that $N$ has  at least $3$ orbits on vertices and
$\Gamma_N$ has diameter at most $2$.
Then $\Gamma$ is a cover of $\Gamma_N$
and either  $\Gamma,\Gamma_N$ are as in Table \ref{table-cover-diam2}
or $\Gamma$ is a non-$(G,4)$-distance-transitive graph of diameter at least $4$ and   $\Gamma_N$ is as in Table \ref{table-cover-diam}.

\end{theo}

The graphs $\Gamma$ and $\Gamma_N$ in Tables \ref{table-cover-diam2} and \ref{table-cover-diam} will be described  in Section 2.
In particular,  `SDC'  denotes  the `standard double cover' (Definition \ref{double-def}).
We remark that in Theorem \ref{girth45-redtheo-1}, if $\Gamma_N$ is strongly regular, then $\Gamma$ and $\Gamma_N$ may have distinct girths. For example,  the Armanios-Wells graph is $(G,3)$-geodesic-transitive of  girth 5 and it is a cover of the  folded $5$-cube $\Box_5$ which has girth 4.
Note that the cycle $C_r$, $r=4$ or 5, is a 2-arc-transitive graph of diameter 2. If $\Gamma$ is a 3-geodesic-transitive cover of $C_r$, then $\Gamma$ has valency 2, and hence $\Gamma$  is a cycle with girth at least 8, and so is 3-arc-transitive.


In Theorem \ref{3gt-redtheo-1}, we could choose  $N$ to be an  intransitive normal subgroup of $G$ which is maximal with respect to having  at least $3$ orbits on the vertex set of $\Gamma$. Then
each  non-trivial normal subgroup $M$ of $G$, properly containing $N$,
has 1 or 2 orbits, and so $M/N$ has 1 or 2 orbits   on the vertices  of $\Gamma_N$.
In other words,   $G/N$ is quasiprimitive or bi-quasiprimitive on the vertex set of $\Gamma_N$.
Hence this theorem suggests that to investigate the family of $(G,3)$-geodesic-transitive graphs that are not $(G,3)$-arc-transitive,  we should
concentrate on the following two problems:

\begin{prob}\label{prob-1}
{\rm
\begin{itemize}
\item[(1)]  determine  $(G,3)$-geodesic-transitive graphs of girth 4 or 5 where $G$ acts quasiprimitively or bi-quasiprimitively
on the vertex set;

\item[(2)] investigate $(G,3)$-geodesic-transitive covers of
the graphs obtained from  (1),
and also investigate  the graphs $\Gamma$ in Table \ref{table-cover-diam}.

\end{itemize}
}
\end{prob}

This paper is organised  as follows.
After this Introduction, we give, in Section 2,
some definitions on groups and graphs that we need and also prove some elementary  lemmas which will be used in the following analysis.
Theorem \ref{3gt-redtheo-1} is proved in   Section 3. It
reduces  the study of $3$-geodesic-transitive graphs of  girth 4 or 5 to the study of $G$-normal covers of complete graphs, of strongly
regular graphs, and of a class of $3$-geodesic-transitive graphs of diameter at least 3 with the same girth.
Then in Section 4, we determine all the $2$-arc-transitive  strongly regular graphs and complete the proof of
Theorem \ref{diam2-theo-1}.
In Section 5, we prove Theorem \ref{girth45-redtheo-1}, that is, we
investigate  3-geodesic-transitive graphs of girth 4 or 5 which are   covers of a 2-arc-transitive graph that is complete or strongly regular.

\section{Preliminaries}

In this section, we give some definitions concerning groups and graphs and also prove some  results which will be used in our
analysis.

\subsection{Graph theoretic notions} All graphs in this paper are finite, simple, connected and undirected. For a graph $\Gamma$, we use $V(\Gamma)$ and $\Aut(\Gamma)$ to denote
its \emph{vertex set}  and \emph{automorphism
group}, respectively. For the group theoretic terminology not defined here we refer the reader to \cite{Cameron-1,DM-1,Wielandt-book}.

In  a graph $\Gamma$, $d_{\Gamma}(u,v)$ denotes the distance between two vertices $u$ and $v$ in $\Gamma$. The \emph{diameter} $\diam(\Gamma)$ of  $\Gamma$ is the maximum distance between $u,v$
for all $u, v \in V(\Gamma)$.
If  $\Gamma$ is $(G,s)$-geodesic-transitive with $G=\Aut(\Gamma)$, then $\Gamma$ is  called \emph{$s$-geodesic-transitive}; and if further
$s=\diam(\Gamma)$, then $\Gamma$ is  called \emph{geodesic-transitive}.
A graph $\Gamma$ is said to be  \emph{$(G,s)$-distance-transitive} ($(G,s)$-DT) if $G\leq \Aut(\Gamma)$ and for each $i\leq s$, $G$  is transitive on
the ordered pairs of vertices at  distance $i$. Moreover,     a $(G,s)$-distance-transitive graph is said to be \emph{$G$-distance-transitive} ($G$-DT) if $s=\diam(\Gamma)$, and
\emph{distance-transitive}  if $G=\Aut(\Gamma)$.
By  definition, every $(G,s)$-geodesic-transitive graph is  $(G,s)$-distance-transitive.

Let  $G$ be a transitive permutation group on a set  $\Omega$.
Let $B$ be a non-empty
subset of $\Omega$. Then $B$ is called a \emph{block} of
$G$ if, for any $g\in G$, either $B^g=B$ or
$B^g\cap B =\emptyset$. The set   $\Omega$ and singleton subsets are trivial blocks.
If $N$ is a normal subgroup of $G$, then each $N$-orbit is a block of $G$ on $\Omega$, and the
set of $N$-orbits    forms a
$G$-invariant partition of $\Omega$.
The group $G$ is said to be \emph{primitive} on $\Omega$, if $G$
 has only trivial blocks.  There is a remarkable classification of finite primitive permutation
groups (8 types), mainly due to M. O'Nan  and L. Scott, see
\cite{LPS-1}.

A subgraph $X$ of $\Gamma$ is an \emph{induced subgraph} if two
vertices of $X$ are adjacent in $X$ if and only if they are adjacent
in $\Gamma$.  When $U\subseteq V(\Gamma)$, we use $[U]$ to denote  the
subgraph of $\Gamma$ induced by $U$.
Let $u$ be a vertex in a graph $\Gamma$ and let $i\leq \diam(\Gamma)$.
Then  $\Gamma_i(u)$ denotes the set of vertices at  distance $i$ from $u$.
In the characterisation of $s$-geodesic-transitive  and $s$-distance-transitive
graphs, the following parameters are important.

\begin{defi}\label{intersectionarray}
{\rm Let $\Gamma$ be an  $s$-distance-transitive graph, let $u\in V(\Gamma)$, and
let $v\in \Gamma_i(u)$, $i\leq s$. Then the number of edges
from $v$ to $\Gamma_{i-1}(u)$, $\Gamma_i(u)$, and $\Gamma_{i+1}(u)$
does not depend on the choice of $v$ and these numbers are
denoted, respectively, by $c_i$, $a_i$, $b_i$.}
\end{defi}

Clearly  $a_i+b_i+c_i$ is equal to the valency of $\Gamma$ whenever the constants are well-defined. Note that for $3$-geodesic-transitive graphs, the constants are always well-defined for $i=1, 2, 3$.
If $\Gamma$ is distance-transitive then the constants are  well-defined for $i=1, 2, \ldots,d$ where $d$
is the diameter of $\Gamma$, and
the sequence $(b_0,\ldots,b_d;c_1,\ldots,c_d)$,  is called the \emph{intersection array} of $\Gamma$.
Some properties of these parameters are given in \cite[Proposition 20.4]{Biggs-1}.

\subsection{The graphs occurring in Theorem  \ref{girth45-redtheo-1}}

The \emph{Hoffman-Singleton graph} $HoS$, is the unique strongly regular graph  with parameters $(50,7,0,1)$, automorphism group $PSU(3,5).\mathbb{Z}_2$ and vertex stabiliser $S_7$.
The \emph{Higman-Sims graph} $HiS$,  is the unique strongly regular graph with parameters $(100,22,0,6)$,  automorphism group $HS.\mathbb{Z}_2$ and vertex stabiliser $M_{22}.\mathbb{Z}_2$.
The \emph{Gewirtz graph}  is the unique strongly regular graph with parameters $(56,10,0,2)$, automorphism group $PSL(3,4).\mathbb{Z}_2^2$ and vertex stabiliser $A_6.\mathbb{Z}_2^2$.

The    \emph{Hamming graph}  $\H(d,2)$  has
vertex set $\Delta^d=\{(x_1,x_2,\ldots,x_d)|x_i\in \Delta\}$, where
$\Delta=\{0,1\}$, and  two vertices $v$ and $v'$ are adjacent if and only if they differ
in exactly one coordinate. The    Hamming graph  $\H(d,2)$ is also called a \emph{$d$-cube};
it  is  3-geodesic-transitive  with girth 4 and diameter $d$ whenever  $d\geq 3$.
The \emph{folded $d$-cube} $\Box_d$, is the graph obtained by identifying antipodal  vertices of $\H(d,2)$, i.e., vertices  at  distance $d$.
Hence  $\H(d,2)$ is the antipodal cover of  $\Box_d$ with antipodal parts of  size 2.
The folded $d$-cube has diameter $\lfloor \frac{d}{2} \rfloor$ and valency $d$, and if $d\geq 4$, then  $\Box_d$ has  girth 4.

 Let $\Gamma$ be the  Hoffman-Singleton graph and let $u\in V(\Gamma)$. Then the induced subgraph $[\Gamma_2(u)]$, denoted by $[HoS]_2$, is a distance-transitive graph with diameter 3, girth 5, see \cite[p.223]{BCN}, and it is  a 3-geodesic-transitive  6-cover of $\K_7$.
Let $\Gamma=HiS$  and let $u\in V(\Gamma)$. Then the induced subgraph $[\Gamma_2(u)]$ is called the \emph{$M_{22}$-graph}.
It is   the unique strongly regular graph with parameters $(77, 16, 0, 4)$,  automorphism group  $M_{22}.\mathbb{Z}_2$ and vertex stabiliser  $\mathbb{Z}_2^4:S_6$,
see \cite[p.369]{BCN}.

A   distance-regular graph with intersection array $(2\mu,2\mu-1,\mu,1;1,\mu,2\mu-1,2\mu)$
is called a \emph{ Hadamard graph} of valency
$2\mu$ (for definition of distance-regular graph  see \cite[page 1]{BCN}, and intersection array is defined in the next subsection).
Hence a Hadamard graph of valency $2\mu$ has $8\mu$ vertices.



The \emph{Armanios-Wells graph} is the unique distance-regular graph with   intersection array $(5,4,1,1;1,1,4,5)$,  automorphism group $\mathbb{Z}_2^{1+4}.A_5$, and vertex stabiliser  $A_5$.
It   is 2-arc-transitive, and also 3-geodesic-transitive.
The \emph{dodecahedron} has both girth  and diameter 5,   automorphism group $A_5\times \mathbb{Z}_2$ and intersection array  $(3,2,1,1,1;1,1,1,2,3)$, see \cite[p.1]{BCN} and \cite{Frucht-1};
it is geodesic-transitive.

A \emph{divisible design} $GD(k,\lambda,n,kn)$ is a triple $(X,\mathcal{P},\mathcal{B})$
where $X$ is a set of $kn$ `points', $\mathcal{P}$ is a partition of $X$ into classes of size $n$, and
$\mathcal{B}$ is a collection of $k$-subsets of $X$ (called `blocks') such that each block meets every class
in precisely 1 point, and any two points of $X$  from different classes are contained in $\lambda$
blocks. A design is called \emph{resolvable} when its set of blocks can be partitioned into parallel classes, that is,
into partitions of the point set. We use $RGD(k,\lambda,n)$ to denote a resolvable divisible design $GD(k,\lambda,n,kn)$
(see \cite[page 439]{BCN}).


Let   $\mathcal{D}=(X,\mathcal{P},  \mathcal{B})$ be a resolvable divisible   design $RGD(r,\lambda,m)$ such that $r=\lambda.m$.
Counting triples $(u,v,B)$, such that $u,v$ are distinct points contained in $B\in\mathcal{B}$, yields $|\mathcal{B}|=m^2\lambda = rm
=|\mathcal{P}|$. Thus each parallel class of blocks contains $m$ blocks, and it follows that $\mathcal{B}$ is a disjoint union of $r$ parallel classes of blocks,
and then, since each point lies in exactly one block from each parallel class,  each point lies in exactly $r$ blocks, that is to say,  $\mathcal{D}$  is a 1-$(rm,r,r)$ design. Let $\mathcal{B}=\{B_1,B_2,\ldots,B_{rm}\}$. The \emph{incidence graph} $\Inc(\mathcal{D})$ of $\mathcal{D}$ is defined as follows:
the  vertex set is $V_1\cup V_2$ where $V_1=X$ and $V_2=\mathcal{B}$, and a vertex  $x$ of $V_1$ is adjacent to a vertex  $B_i$ of $V_2$
if and only if $x\in B_i$.

\begin{lemma}\label{res-divisible-1}
Let   $\mathcal{D}=(X,\mathcal{P},  \mathcal{B})$ be a resolvable divisible   design $RGD(r,\lambda,m)$ with $r=\lambda.m$,
and define $\Inc(\mathcal{D})$ as above. Suppose moreover that  any two  blocks from different parallel classes contain exactly
$\lambda$ common points. Then   $\Inc(\mathcal{D})$ is a  distance-regular bipartite antipodal cover of $\K_{r,r}$ with  diameter $4$ and antipodal block  size $m$, written $\Inc(\mathcal{D}) = m\K_{r,r}$.
\end{lemma}

\proof By definition  $\Inc(\mathcal{D})$ is a bipartite  graph with bipartition $\{ V_1, V_2\}$.
Moreover, $|V_1|=|V_2|=rm$.
Now  $\mathcal{P}$ is a partition of $V_1$ into $r$ classes of size $m$,  say $\mathcal{P}=\{P_1,P_2,\ldots,P_r\}$, where each $|P_i|=m$.
As any two points of $V_1$  from different classes are contained in $\lambda$
blocks, it follows that vertices of $V_1$ from distinct classes are at distance 2 in  $\Inc(\mathcal{D})$.

Also  $\mathcal{B}$ is a disjoint union of $r$ parallel classes of blocks.
Let $\mathcal{C}=\{C_1,C_2,\ldots,C_r\}$ be the set of parallel classes, where each $C_i$ contains $m$ blocks of $\mathcal{B}$.
By assumption, any two blocks from different parallel classes contain exactly $\lambda$ common points, and hence
vertices of $V_2$ from distinct parallel classes are at distance 2 in  $\Inc(\mathcal{D})$.
Further,  for each $i$, each vertex of $V_1$ is adjacent to exactly one vertex of $C_i$,
and each vertex of $V_2$ is adjacent to exactly one vertex of $P_i$.

Let $x\in V_1$ and $C\in\mathcal{C}$. Then $x$ lies in a unique block $B\in C$. Let $B''\in C\setminus\{B\}$, and $C'\in\mathcal{C}\setminus\{C\}$,
so  $x$ also lies in a unique block $B'\in C'$. By assumption $|B'\cap B''|=\lambda$. Let $y\in B'\cap B''$. Then $(x, B', y, B'')$ is a path of length 3 in
 $\Inc(\mathcal{D})$, and hence the distance between $x$ and $B''$ in  $\Inc(\mathcal{D})$ is $3$.  It follows that $x$ has distance 1 or 3 from
 each vertex of $V_2$. If $x\in P\in \mathcal{P}$, then as mentioned above, $x$ is at distance 2 from each point
 $y\in V_1\setminus P$. If $x'\in P\setminus\{x\}$, then $x'$ also has distance $2$ from $y$, and hence $x, x'$ are at distance $4$ in
  $\Inc(\mathcal{D})$. By a similar argument, each vertex of $V_2$  is at distance 1 or 3 from vertices of $V_1$, and at distance 2 or 4 from vertices of $V_2$.
Thus  $\Inc(\mathcal{D})$ has diameter 4, and  two vertices of  $\Inc(\mathcal{D})$ are at distance 4 if  and only if they are distinct blocks in the same parallel class, or distinct points in the same class $P$.
Hence  the classes of $V_1$ and the parallel  classes of $V_2$ are antipodal blocks of  $\Inc(\mathcal{D})$.
It is straightforward, using the properties of the design $\mathcal{D}$,  to show that  $\Inc(\mathcal{D})$  is distance-regular with intersection array
$(r, r-1, r-\lambda, 1; 1, \lambda, r-1, r)$.

Finally, we examine the antipodal quotient graph, that is, the graph $\Gamma'$ with vertex set $\mathcal{P}\cup \mathcal{C}$, such that
a vertex $P_i\in \mathcal{P}$ is adjacent to a vertex $C_j\in \mathcal{C}$
if and only if some vertex of $P_i$ is adjacent to some vertex of  $C_j$ in  $\Inc(\mathcal{D})$.
Since, in fact, each vertex $x\in P_i$ is contained in a unique block in $C_j$, and vice versa, it follows that $\Gamma'\cong \K_{r,r}$.
Moreover, since $|P_i|=|C_j|$, we see that the edges of  $\Inc(\mathcal{D})$ joining vertices of $P_i$ and $C_j$ form a perfect matching
between these two sets. Thus  $\Inc(\mathcal{D})$ is a cover of $\Gamma'$.
 \qed

\subsection{Some results}

We will use the following remark frequently throughout the paper. It is  especially useful for studying  small valency $s$-geodesic-transitive graphs.
Recall Definition~\ref{intersectionarray}.

\begin{rem}\label{3gt-dt-1}
{\rm
(1) Let $\Gamma$ be a $2$-geodesic-transitive graph with $b_2\leq 1$.
If   $b_2=0$, then as $\Gamma$ is $2$-geodesic-transitive,
$\Gamma$ has diameter 2 and so  is geodesic-transitive.
Suppose that $b_2=1$, so in particular, $d:=\diam(\Gamma)\geq 3$. Let $(u_0,\ldots,u_d)$ be a $d$-geodesic. Then for each $j\leq d-3$, it follows from  the $2$-geodesic-transitivity of $\Gamma$ that  $|\Gamma_3(u_j)\cap \Gamma(u_{j+2})|=1$.
Note that, $\Gamma_{j+3}(u_0)\cap \Gamma(u_{j+2})\subseteq \Gamma_3(u_{j})\cap \Gamma(u_{j+2})$,
and so
$|\Gamma_{j+3}(u_0)\cap \Gamma(u_{j+2})|=1$. Hence the 2-geodesic $(u_0,u_1,u_2)$ has a unique continuation to an $r$-geodesic in $\Gamma$ for all $r$ such that $3\leq r\leq d$.
Since $\Gamma$ is $2$-geodesic-transitive,  each 2-geodesic of $\Gamma$ has a unique continuation to an $r$-geodesic for  $3\leq r\leq d$.
Thus  $\Gamma$ is  geodesic-transitive, and hence also distance-transitive.

(2) Let $\Gamma$ be a $3$-geodesic-transitive graph with $b_3\leq 1$.
If   $b_3=0$, then since $\Gamma$ is $3$-geodesic-transitive,
$\Gamma$ has diameter 3 and so  is geodesic-transitive.
Suppose that $b_3=1$. So in particular, $d:=\diam(\Gamma)\geq 4$.
Let $(u_0,\ldots,u_d)$ be a $d$-geodesic. Then for each $j\leq d-4$, it follows from the $3$-geodesic-transitivity of $\Gamma$ that   $|\Gamma_4(u_j)\cap \Gamma(u_{j+3})|=1$.
Note that, $\Gamma_{j+4}(u_0)\cap \Gamma(u_{j+3})\subseteq \Gamma_4(u_{j})\cap \Gamma(u_{j+3})$, and so
$|\Gamma_{j+4}(u_0)\cap \Gamma(u_{j+3})|=1$.
Hence the 3-geodesic $(u_0,u_1,u_2,u_3)$ has a unique continuation to an $r$-geodesic in $\Gamma$ for all $r$ such that $4\leq r\leq d$.
Since $\Gamma$ is $3$-geodesic-transitive,  each 3-geodesic of $\Gamma$ has a unique continuation to an $r$-geodesic for $4\leq r\leq d$.
Thus  $\Gamma$ is geodesic-transitive, and hence also distance-transitive.

 }
\end{rem}

\begin{defi}\label{double-def}
{\rm  For   be a graph $\Gamma$ with vertex set $V(\Gamma)$ and arc set $A(\Gamma)$, let
$\overline{\Gamma}$ be the   graph with   vertex set $V(\Gamma)\times \{1,2\}$, such that two vertices $(x,1)$ and $(y,2)$  adjacent if and only if $(x,y)\in A(\Gamma)$. Then  $\overline{\Gamma}$ is called the  \emph{standard double cover} (SDC) of $\Gamma$; $\overline{\Gamma}$ is bipartite with bipartite halves $V(\Gamma)\times \{i\}$ for  $i=1,2$.


  }
\end{defi}

Some of our examples are standard double covers of distance-transitive graphs, and we use the following observation to identify them.

\begin{lemma}\label{sdc}
Suppose that $\Gamma$  is a finite distance-transitive graph of odd diameter which is both bipartite and antipodal with antipodal blocks of size $2$. Let $\Sigma$ denote the antipodal quotient of $\Gamma$. Then $\Gamma$ is isomorphic to the standard double cover of $\Sigma$.
\end{lemma}

\proof
Let $d=\diam(\Gamma)$, let $\Delta_1, \Delta_2$ be the parts of the bipartition of $\Gamma$, and identify $V(\Sigma)$ with the set of antipodal blocks. Let $B=\{u,v\}\in V(\Sigma)$. Then $\Gamma_d(u)=
\{v\}$, and since $d$ is odd we have, say, $u\in\Delta_1, v\in\Delta_2$. Let $C\in V(\Sigma)$ be adjacent to $B$ in $\Sigma$.
Then $C=\{x,y\}$, and $u$ is adjacent to exactly one vertex of $C$, say $x$. Thus $x\in\Gamma(u)\subseteq \Delta_2$ and $y\in\Gamma(v)=\Gamma_{d-1}(u)\subseteq \Delta_1$.

Define a map $\varphi: V(\Gamma)\mapsto V(\Sigma)\times\{1,2\}$ by $\varphi(u)=(B, i)$ where $B$ is the antipodal block containing $u$ and $\Delta_i$ is the part of the bipartition containing $u$. Then $\varphi$ is a well defined bijection. We claim that $\varphi$ is a graph isomorphism from $\Gamma$ to $\overline{\Sigma}$, the standard double cover of $\Sigma$: for if $\{u,x\}$ is an edge of $\Gamma$ and $\varphi(u)=(B,i), \varphi(x)=(C,j)$ then $i\ne j$ (as $u,x$ lie in different parts of the bipartition) and $\{B, C\}$ is an edge of $\Sigma$ (by definition of $\Sigma$). Conversely if
$\varphi(u)=(B,i), \varphi(x)=(C,j)$ form an edge of $\overline{\Sigma}$, then $i\ne j$ and $\{B, C\}$ is an edge of $\Sigma$, by the definition of  $\overline{\Sigma}$. As $\Gamma$ covers $\Sigma$, there are exactly two edges of $\Gamma$ joining vertices of $B$ and $C$. Now $u$ is the
unique vertex in $B\cap\Delta_i$, and since $\Gamma$ is bipartite, the edge involving $u$  must join $u$ to the unique vertex in $C\cap\Delta_j$, namely $x$. This proves the claim and hence $\Gamma\cong\overline{\Sigma}$.
\qed

If $G\leq \Aut(\Gamma)$, then $G$ also acts as a group of automorphisms of the standard double cover $\overline{\Gamma}$ with
the action $g: (x, i)\mapsto  (x^g, i)$. If $G$ is vertex-transitive on $\Gamma$, then $G$ has two orbits
on the set of vertices of $\overline{\Gamma}$, namely the sets $V(\Gamma)\times \{i\}$ for $i=1,2$.
Furthermore, $G_v = G_{(v,i)}$ for each $i = 1, 2$ and $v\in V(\Gamma)$,
the action of $G_v$ on $\Gamma(v)$ is equivalent to the action of
$G_{(v,i)}$ on $\overline{\Gamma}((v,i))$, and  if $\Gamma$ is $(G,2)$-arc-transitive, then $\overline{\Gamma}$ is locally $(G,2)$-arc-transitive.
Define
\begin{align*}
\tau:&V(\overline{\Gamma})\mapsto V(\overline{\Gamma}),\\
&(v,i)\mapsto (v,3-i).
\end{align*}

Then $\tau$ is a graph automorphism of $\overline{\Gamma}$ of  order 2. Further, for any vertex $(v,i)$, we have $(v,i)^{g\tau}=(v^g,i)^\tau=(v^g,3-i)=(v,3-i)^g=(v,i)^{\tau g}$.
Hence  $g\tau=\tau g$ for every $g\in G$. Let $\overline{G}=G\times \langle \tau \rangle$. Then $\overline{G}\leq \Aut(\overline{\Gamma})$.
Moreover, $\Gamma$ is $(G,2)$-arc-transitive, then $\overline{\Gamma}$ is $(\overline{G},2)$-arc-transitive.

Now $(x, 1)$ is adjacent to $ (y, 2)$  if and only if $(y, 1) $ is adjacent to $ (x, 2)$. If $\Gamma$ is  connected, then for distinct $x,y\in V(\Gamma)$ there exists a path $P$ in $\Gamma$
from  $x$ to $y$. The path $P$ lifts to a path in $\overline{\Gamma}$ from  $(x, 1)$ to $(y, 1)$ if $P$ has
even length, and to one between $(x, 1)$ and $(y, 2)$ if $P$ has odd length. In particular, there is a
path from  $(y, 1)$ to $(y, 2)$ if and only if $y$ lies  in a cycle  in $\Gamma$ of odd length. Thus  $\overline{\Gamma}$ is connected if and only if $\Gamma$ is connected and  contains an odd cycle,
that is, if and only if $\Gamma$ is connected and not bipartite.

We mention a few more properties of these standard double covers. When we say that a (non-complete)
graph $\Gamma$  `\emph{has $c_2=c$}' we mean that $|\Gamma(u)\cap \Gamma(w)|=c$ whenever
$d_\Gamma(u,w)=2$.

\begin{lemma}\label{rect-cover-1}
Let $\Gamma$ be a    non-bipartite graph of girth $4$ with vertex-transitive group $G\leq\Aut(\Gamma)$,  and let $\overline{\Gamma}$ be its standard double cover with group $\overline{G}$ as above.
Then the following hold.

 \begin{itemize}
\item[(1)] For an integer $c\geq2$, $\Gamma$ has $c_2=c$  if and only if   $\overline{\Gamma}$ has   $c_2=c$.

\item[(2)] For a positive integer $s\leq\diam(\Gamma)$, $\Gamma$ is $(G,s)$-distance-transitive   if and only if   $\overline{\Gamma}$ is $(\overline{G},s)$-distance-transitive.

\item[(3)] If $a_2=0$ for $\Gamma$, then $\Gamma$ is $(G,3)$-geodesic-transitive   if and only if   $\overline{\Gamma}$ is $(\overline{G},3)$-geodesic-transitive.

\end{itemize}

\end{lemma}
\proof (1) Suppose that $\Gamma$ has  $c_2=c$.  Let $d_\Gamma(u,w)=2$ and, say, $\Gamma(u)\cap \Gamma(w)=\{v_1,\dots, v_c\}$. Then for each $i$,
$[(u,1),(v_i,2),(w,1)]$ is a 2-arc of $\overline{\Gamma}$, so $\{(v_1,2),\dots,(v_c,2)\} $ $ \subseteq \overline{\Gamma}((u,1))\cap \overline{\Gamma}((w,1))$. Equality holds, since if
$ \overline{\Gamma}((u,1))\cap \overline{\Gamma}((w,1))$ contains a vertex $(u,2)$, then $u\in \Gamma(u)\cap \Gamma(w)$. Similarly $\{(v_1,1),\dots,(v_c,1)\}= \overline{\Gamma}((u,2))\cap \overline{\Gamma}((w,2))$. Thus $\overline{\Gamma}$ has parameter $c_2=c$.

Conversely, suppose that $\overline{\Gamma}$ has parameter $c_2=c$.  A distance-two pair is of the form $(u,i), (w,i)$ for vertices $u,w\in V(\Gamma)$ and $i\in\{1,2\}$. Then
$\overline{\Gamma}((u,i))\cap \overline{\Gamma}((w,i))= \{(v_1,3-i),\dots,(v_c,3-i)\}$ and a similar argument to the above shows that $\Gamma(u)\cap \Gamma(w)=\{v_1,\dots,v_c\}$, and so  $\Gamma$ has parameter $c_2=c$.

(2) Since $G, \overline{G}$ are vertex-transitive on $\Gamma, \overline{\Gamma}$, respectively,
and since, for $u\in V(\Gamma)$ and $j\in\{1,2\}$ we have $G_u = \overline{G}_{(u,j)}$, it is sufficient to
examine the orbits of $G_u$.  For $i\leq s$, $\overline{\Gamma}_i((u,j))= \{ (v,j')\mid v\in\Gamma_i(u)\}$, where $j'=j$ if $i$ is even and $j'=3-j$ if $i$ is odd. The actions of $G_u$ on
$\Gamma_i(u)$ and $ \overline{\Gamma}_i((u,j))$ are therefore equivalent for each $i$, and part (2) follows.

(3) Suppose that $a_2=0$ for $\Gamma$. Assume that $\Gamma$ is $(G,3)$-geodesic-transitive.
Then since $\Gamma$ has  girth 4, $G_u$ is 2-transitive on $\Gamma(u)$, and so
$G_{(u,1)}$ is 2-transitive on $\overline{\Gamma}((u,1))$, hence $\overline{\Gamma}$ is
$(\overline{G},2)$-arc-transitive. Let $[(u,1),(v_1,2),(v_2,1),$ $(u_4,2)]$ and $[(u,1),(v_1,2),(v_2,1),(u_5,2)]$ be two 3-geodesics of $\overline{\Gamma}$. Then
$(u,v_1,v_2,u_4)$ and $(u,v_1,v_2,$ $u_5)$ are two 3-geodesics of $\Gamma$, as  $a_2=0$ for $\Gamma$. Thus some element of $G$ maps
$(u,v_1,v_2,u_4)$ to $(u,v_1,v_2,u_5)$, and so this element maps
$[(u,1),(v_1,2),(v_2,1),(u_4,2)]$ to $[(u,1),(v_1,2),(v_2,1),(u_5,2)]$ in $\overline{\Gamma}$. Hence $\overline{\Gamma}$ is $(G,3)$-geodesic-transitive.

Conversely, suppose that $\overline{\Gamma}$ is $(\overline{G},3)$-geodesic-transitive.
Then as  $\overline{\Gamma}$ is bipartite and
$(\overline{G},2)$-arc-transitive, $\overline{G}_{(u,1)}$ is 2-transitive on $\overline{\Gamma}((u,1))$, and so $G_u$ is 2-transitive on $\Gamma(u)$,
hence  $\Gamma$ is $(G,2)$-arc-transitive. Let
$(u,v_1,v_2,u_4)$ and $(u,v_1,v_2,u_5)$ be two 3-geodesics of $\Gamma$.
Then as $a_2=0$ for $\Gamma$,  $[(u,1),(v_1,2),$ $(v_2,1),(u_4,2)]$ and $[(u,1),(v_1,2),$ $(v_2,1),(u_5,2)]$ are two 3-geodesics of $\overline{\Gamma}$.  Thus
some element of $\overline{G}$ maps
$[(u,1),(v_1,2),$ $(v_2,1),(u_4,2)]$ to $[(u,1),(v_1,2),(v_2,1),(u_5,2)]$,
this element fixes $V(\Gamma)\times \{1\}$ and so lies in $G$,
and hence
this element induces an element of $G$ that maps
$(u,v_1,$ $v_2,u_4)$ to $(u,v_1,v_2,u_5)$. Thus $\Gamma$ is $(G,3)$-geodesic-transitive.
\qed

\section{Reduction result}

In this section, we study normal quotients of  $3$-geodesic-transitive  graphs of  girth 4 or 5.
We will need the following  result from \cite[Lemma 5.3]{DGLP-locdt-2012}.

\begin{lemma}\label{dt-quotient12}
Let $\Gamma$ be a connected locally $(G, s)$-distance-transitive
graph with $s\geq 2$. Let $1\neq N \lhd G$ be intransitive on
$V(\Gamma)$, and let $\mathcal{B}$ be the set of $N$-orbits on
$V(\Gamma)$. Then one of the following holds:

 \begin{itemize}
\item[(i)] $|\mathcal{B}| = 2$.

\item[(ii)] $\Gamma$ is bipartite, $\Gamma_N\cong \K_{1,r}$ with
$r\geq 2$ and $G$ is intransitive on $V(\Gamma)$.

\item[(iii)] $s=2$, $\Gamma\cong \K_{m[b]}$,   $\Gamma_N \cong \K_{m}$
with $m\geq 3$ and $b\geq 2$.

\item[(iv)] $N$ is semiregular on $V(\Gamma)$,  $\Gamma$ is a cover
of $\Gamma_N$, $|V(\Gamma_N)|<|V(\Gamma)|$ and $\Gamma_N$ is
locally $(G/N,s')$-distance-transitive where $s'=\min\{s,\diam(\Gamma_N)\}$.

\end{itemize}
\end{lemma}

We derive from Lemma~\ref{dt-quotient12} the following result for $s$-geodesic-transitive graphs.

\begin{lemma}\label{2gt-covergt-1}
Let $\Gamma$ be a connected $(G,s)$-geodesic-transitive graph  where $s\geq
2$.  Let $1\neq N \lhd G$ be intransitive on $V(\Gamma)$.
Suppose that $\Gamma\ncong \K_{m[b]}$ for any   $m\geq 3$ and $b\geq 2$. Then either

 \begin{itemize}
\item[(i)] $N$ has $2$ orbits on
$V(\Gamma)$ and $\Gamma$ is bipartite; or

\item[(ii)]  $N$  has at least $3$ orbits on
$V(\Gamma)$, $N$ is semiregular on $V(\Gamma)$, $\Gamma$ is a cover of
$\Gamma_N$ and  $\Gamma_N$ is $(G/N,s')$-geodesic-transitive where
$s'=\min\{s,\diam(\Gamma_N)\}$.

\end{itemize}

\end{lemma}

\proof By assumption $N$ is not transitive on $V(\Gamma)$. If $N$ has exactly $2$ orbits on
$V(\Gamma)$,  say  $C_0$ and $C_1$,  then as $\Gamma$ is connected and $G$-arc-transitive,
each $C_i$  contains no edges of $\Gamma$, and so
$\Gamma$ is a bipartite graph with  $C_0,C_1$ being the two bipartite halves, part (i) holds. Thus, we may  suppose that
$N$  has at least $3$ orbits on
$V(\Gamma)$.   Since $\Gamma$ is  $(G,s)$-geodesic-transitive   where $s\geq
2$, $\Gamma$ is also $(G,s)$-distance-transitive.
Since $\Gamma\ncong \K_{m[b]}$ for any   $m\geq 3$ and $b\geq 2$, it follows that part (iv) of
Lemma \ref{dt-quotient12} holds. To prove that part (ii) is valid it remains to prove that
$\Gamma_N$ is locally $(G/N,s')$-geodesic-transitive, where $s'=\min\{s,\diam(\Gamma_N)\}$.

Let  $(B_0,B_1,B_2,\ldots,B_{t})$ and
$(C_0,C_1,C_2,\ldots,C_{t})$ be  $t$-geodesics of $\Gamma_N$
where $t\leq s'$. Since $\Gamma$ is a
cover of $\Gamma_N$, there exist $x_i\in B_i$ and $y_i\in C_i$ such
that $(x_0,x_1,x_2,\ldots,x_{t})$ and $(y_0,y_1,y_2,\ldots,y_{t})$
are $t$-geodesics of $\Gamma$. As $t\leq s'\leq s$ and $\Gamma$
is $(G,s)$-geodesic-transitive,  there exists $g\in G$ such that
$(x_0,x_1,x_2,\ldots,x_{t})^g=(y_0,y_1,y_2,\ldots,y_{t})$, and hence
$g$ maps
$(B_0,B_1,B_2,\ldots,B_{t})$ to $(C_0,C_1,$ $C_2,\ldots,C_{t})$. Thus
$\Gamma_N$ is $(G/N,s')$-geodesic-transitive, and (ii) holds. \qed

\begin{lemma}\label{girth45-lem-1}
Let $\Gamma$ be a   $(G,2)$-geodesic-transitive  graph of girth at least $4$.
Let $N$ be an intransitive normal subgroup of $G$ with  at least $3$ orbits on
$V(\Gamma)$. Then  $\Gamma_N$ is a complete graph if and only if  $\Gamma_N$ has girth $3$.

\end{lemma}

\proof  Note that $\Gamma\ncong \K_{m[b]}$ for any   $m\geq 3$ and $b\geq 2$, since
the   girth of $\Gamma$ is  at least $4$.
If $\Gamma_N$ is a complete graph, then since  $N$  has at least $3$ orbits on
$V(\Gamma)$,  it follows that $\Gamma_N$ has girth $3$. Conversely, suppose that   $\Gamma_N$ has girth $3$.
Then by  Lemma \ref{2gt-covergt-1}, the graph   $\Gamma$ is a cover of $\Gamma_N$.
Let $(B_1,B_2,B_3)$ be a triangle of $\Gamma_N$.
Then   there exist $b_i\in B_i$ such that   $(b_1,b_2,b_3)$ is a 2-arc of $\Gamma$. Further, $(b_1,b_2,b_3)$ is a 2-geodesic of $\Gamma$ (as $\Gamma$ has girth at least $4$).
Suppose that $\Gamma_N$ is not a complete graph. Then $\Gamma_N$ has a 2-geodesic $(C_1,C_2,C_3)$.  Since  $\Gamma$  covers $\Gamma_N$,  we can find  $c_i\in C_i$ such that  $(c_1,c_2,c_3)$  is a 2-geodesic of $\Gamma$.
As $\Gamma$ is  $(G,2)$-geodesic-transitive,  there exists $g\in G$
such that $(b_1,b_2,b_3)^g=(c_1,c_2,c_3)$, and so  $(B_1,B_2,B_3)^g=(C_1,C_2,C_3)$, which is impossible.
Therefore $\Gamma_N$ is a complete graph. \qed

\medskip
Now we  prove Theorem \ref{3gt-redtheo-1} and
describe the various possibilities for the girth and diameter of $\Gamma_N$.

\medskip
\noindent {\bf Proof of Theorem \ref{3gt-redtheo-1}.}
 Since $\Gamma$ has girth 4 or 5, we have $\Gamma\ncong \K_{m[b]}$ for any   $m\geq 3$ and $b\geq 2$.
Thus by   Lemma \ref{2gt-covergt-1},   $\Gamma$ is a cover of $\Gamma_N$ and $\Gamma_N$ is $(G,s')$-geodesic-transitive where
$s'=\min\{3,\diam(\Gamma_N)\}$.
If $\Gamma_N$ is  a complete graph, then (1) holds. Assume now  that $\Gamma_N$ is not a complete graph.
Then by  Lemma \ref{girth45-lem-1},   $\Gamma_N$ has girth at least 4.
Moreover, since $\Gamma_N$ is  covered by  $\Gamma$, it follows that the  girth of $\Gamma_N$ is at most the girth of $\Gamma$,
and hence   $\Gamma_N$ has girth 4 or 5.
If $\Gamma_N$ has  diameter $2$,  then  (2) follows.

Now suppose that $\Gamma_N$  has diameter at least $3$. If  $\Gamma$ has girth  4, then by the previous paragraph  $\Gamma_N$ has girth  4. Suppose that $\Gamma$ has girth  5.
Assume  that the  girth of $\Gamma_N$ is  4 and
let $(B_1,B_2,B_3,B_4)$ be a 4-cycle. Then since  $\Gamma$ is a cover of $\Gamma_N$,  there exist $b_i\in B_i$ and  $b'\in B_1$ such that $(b_1,b_2,b_3,b_4,b')$
is a 4-arc of $\Gamma$. Since $\Gamma$ has girth 5, $b_1\not\in \{b_4, b'\}$ and $(b_1,b_2,b_3)$
is a 2-geodesic.  Furthermore, $b_4\in \Gamma_2(b_1) \cup \Gamma_3(b_1)$.
If  $b_4\in \Gamma_2(b_1)$, then there exists $v\in \Gamma(b_1)$ such that $(b_1,v,b_4)$ is a 2-geodesic.
Let $B$ be the $N$-orbit containing $v$.
Then $(B_1,B,B_4)$ is a triangle of $\Gamma_N$,  contradicting the assumption
that $\Gamma_N$ has girth 4.
Thus $b_4\in \Gamma_3(b_1)$, and so  $(b_1,b_2,b_3,b_4)$
is a 3-geodesic of $\Gamma$.
Let $(C_1,C_2,C_3,C_4)$ be  a 3-geodesic of $\Gamma_N$. Then there exist $c_i\in C_i$ such that $(c_1,c_2,c_3,c_4)$
is a 3-geodesic, and  there exists $g\in G$ such that   $(b_1,b_2,b_3,b_4)^g=(c_1,c_2,c_3,c_4)$,
and hence  $(B_1,B_2,B_3,B_4)^g=(C_1,C_2,C_3,C_4)$.  This is impossible, as $(B_1,B_4)$
is an arc but $C_1,C_4$ are at distance 3 in $\Gamma_N$. Therefore,  $\Gamma_N$ has girth 5, and (3) holds.
 \qed

\section{Finite $2$-arc-transitive strongly regular graphs}

The normal quotient graphs in Theorem \ref{3gt-redtheo-1} (2) are
$2$-arc-transitive strongly regular  graphs. In order to  classify $(G,3)$-geodesic-transitive  graphs of girth $4$ or $5$,
we need to know all possibilities for these normal quotients explicitly, and determining them is the aim of this section.
Note that every $2$-arc-transitive strongly regular graph $\Gamma$  has girth 4 or 5.
Also, for each $u\in V(\Gamma)$ and distinct $v, w\in\Gamma(u)$, the triple $(v,u,w)$ is a $2$-arc, so if $\Gamma$ is $(G,2)$-arc-transitive then $G_u$ is $2$-transitive on $\Gamma(u)$. We frequently uese this fact in our proofs. First we gather results from the literature to
determine all the  girth 5  examples.

\begin{lemma}\label{lemma-2atdg5}
Let $\Gamma$ be a $2$-arc-transitive strongly regular  graph  of girth $5$.  Then    $\Gamma$ is $C_5$,
the Petersen graph,  or the Hoffman-Singleton graph.
\end{lemma}
\proof Since $\Gamma$ is a $2$-arc-transitive   graph of diameter $2$ and girth $5$, it follows from \cite[Theorem 6.7.1]{BCN} that  $\Gamma$ has valency 2, 3, 7 or 57.
By \cite{Aschb-1971} and \cite[p.207, Remark (i)]{BCN},
the   valency of $\Gamma$ is not  57, and so $\Gamma$ has valency 2, 3 or 7.
Moreover, by \cite[p.207, Remark (i)]{BCN} or \cite[p.206]{GR}, if $\Gamma$ has valency 2, then $\Gamma$ is $C_5$;
if $\Gamma$ has valency 3, then $\Gamma$ is the Petersen graph; and if $\Gamma$ has valency 7, then $\Gamma$ is the Hoffman-Singleton graph. \qed

\begin{lemma}\label{lemma-2atdegree}
Let $\Gamma$ be a $2$-arc-transitive  strongly regular graph.  Then  either  $\Gamma$ is a complete bipartite graph, or $\Gamma\cong C_5$, or $|\Gamma(u)|<|\Gamma_2(u)|$ for each vertex $u$.
\end{lemma}

\proof Let $k=|\Gamma(u)|$, and note that $k$ is independent of the choice of $u$. Suppose that $k\geq |\Gamma_2(u)|$.
Since $\Gamma$ is  $2$-arc-transitive and not a complete  graph, it follows that $\Gamma$ has girth at least 4, and
so there are $k(k-1)$ edges between $\Gamma(u)$
and $\Gamma_2(u)$. Hence $k(k-1)=c_2\times |\Gamma_2(u)|\leq c_2k$, so
$c_2=k-1$ or $k$.
If $ c_2=k$, then $\Gamma$ is complete bipartite.
Suppose that $ c_2= k-1$. Then $|\Gamma_2(u)|=k$, and since $\Gamma$ has diameter $2$, we obtain $a_2=1$.
Thus for $v\in\Gamma_2(u)$ and $X={\rm Aut}(\Gamma)$,  $X_{uv}$  fixes $\Gamma(u)\setminus \Gamma(v) = \{x\}$, say, and $\Gamma_2(u)\cap\Gamma(v) = \{w\}$, say. Thus $X_{uv}=X_{ux}$ (as they have the same order), and since $X_u$ is $2$-transitive on $\Gamma(u)$, it follows that
$X_{uv}$ is transitive on $\Gamma(u)\cap \Gamma(v)$. Similarly, $X_{uv}=X_{uw}$ and $X_{uv}$ is transitive on $\Gamma(u)\cap \Gamma(w)$. If the valency $k$ is greater than 2 then we must have $\Gamma(u)\cap \Gamma(v)=\Gamma(u)\cap \Gamma(w)$ and $\Gamma$ contains a 3-cycle, which is a contradiction. Hence $k=2$ and $c_2=1$, and therefore $\Gamma\cong C_5$. \qed

The socle of a finite 2-transitive permutation group is either elementary abelian or  a nonregular nonabelian simple group, see \cite[Theorem 4.1B]{DM-1}. Moreover, in the latter case, the socle is primitive, see \cite[p.244]{DM-1}.

\medskip
Now we prove Theorem \ref{diam2-theo-1} to  determine the class of $2$-arc-transitive strongly regular  graphs.


\medskip

\noindent {\bf Proof of Theorem \ref{diam2-theo-1}.}
Since $\Gamma$ is a $2$-arc-transitive   graph of diameter $2$, $\Gamma$ has girth 4 or 5.
If $\Gamma$ has girth 5, then by Lemma \ref{lemma-2atdg5}, case (2) holds. Assume now that
$\Gamma$ has girth 4.
Let  $X:=\Aut(\Gamma)$ and let $u\in V(\Gamma)$.
The     $2$-arc-transitivity  of $\Gamma$ implies that   the stabiliser $X_u$ is transitive on  both $\Gamma(u)$ and   $\Gamma_2(u)$.

\medskip\noindent
\textit{Case 1.  $X$ imprimitive:}
Suppose that $X$ is not primitive on $V(\Gamma)$. Then $X$ has a non-trivial block on $V(\Gamma)$ containing $u$, say $\Delta$. Since
$\Gamma$ is arc-transitive and connected, $\Delta$ does not contain any edge of $\Gamma$.  Thus  $ \Delta\subseteq  \{u\}\cup \Gamma_2(u)$, as $\Gamma$
has  diameter $2$. Since $X_{u}$ fixes the block $\Delta$ setwise
and acts transitively on $\Gamma_2(u)$, it follows that $\Delta =
\{u\}\cup \Gamma_2(u)$.
As $\Gamma$ is  $2$-distance-transitive, it follows from
\cite[Lemma 5.2]{DGLP-locdt-2012} that either $\Gamma\cong \K_{m[b]}$ for some $m\geq 3,b\geq 2$, or $V(\Gamma)\setminus \Delta=\Gamma(u)$
is the only $X$-image of $\Delta$ different from $\Delta$.
However the graph $\Gamma\cong \K_{m[b]}$ has girth 3 (since $m\geq 3$), and so we
conclude that there are just two blocks of imprimitivity, namely $\Delta$ and $\Gamma(u)$, and hence $\Gamma\cong \K_{m,m}$
for some $m\geq 2$ as in part (1).

\medskip\noindent
\textit{Case 2.  $X$ primitive:}
In the remainder, we suppose that  $X$ acts primitively on $V(\Gamma)$, and as $X_u$ is transitive on $\Gamma(u)$ and   $\Gamma_2(u)$,
$X$ is a primitive rank $3$ group. By  $2$-arc-transitivity, we know that  $X_u$ is 2-transitive on   $\Gamma(u)$. It then follows from  \cite[Theorem A]{Praeger-1988-2trans} that $X$ is  primitive on $V(\Gamma)$ of type either
affine or almost simple.
Note that   the 2-transitive group $X_u^{\Gamma(u)}$ induced by $X_u$ on   $\Gamma(u)$ is also of type  either affine or almost simple.


Let the valency of $\Gamma$ be  $n$. Then since $\Gamma$ is  a $2$-arc-transitive   graph of diameter $2$,
it follows that
$$|V(\Gamma)|\leq 1+n+n(n-1)=n^2+1.$$

\medskip\noindent
\textit{Case 2A.  $X$ affine:}
If  $X$ is an affine group, then $\Gamma$ is among the graphs listed in \cite[Table 1]{IP-1} where the Column 3 entry is `p'.   Thus $\Gamma$ is one of the following graphs:
\begin{center}
$\Box_n$, $P_m(a)$ where $m\geq 3$ and $a=1$ or 2, $\Gamma(C_{23})$ or
$\Gamma(C_{22})$.
\end{center}

The graph
$\Gamma(C_{23})$ has $2^{11}$ vertices of  valency 23, however
$23^2+1<2^{11}$,  so $\Gamma$ is not $\Gamma(C_{23})$; $\Gamma(C_{22})$ has $2^{10}$ vertices of  valency 22, and
$22^2+1<2^{10}$, and so  $\Gamma(C_{22})$ is also not a candidate.
The graph
$P_m(a)$ has $2^{m^a}$ vertices of  valency $n=\frac{2^{am}-1}{2^a-1}$.
By \cite[(1.3)]{IP-1}, if $a=1$, then
$P_m(a)=\K_{2^m}$ has diameter 1, which is not the case; if $a=2$ and $m\geq 3$, then an easy  inductive argument yields $n^2+1<2^{m^a}$, and so
$\Gamma$ is not $P_m(a)$.
Finally, for the folded $n$-cube  $\Box_n$, only $\Box_5$  has diameter 2, and $\Box_5$ is an example, as in case (1).

\medskip\noindent
\textit{Case 2B.  $X$ almost simple:}
From now on we suppose that $X$ is an almost simple group with nonabelian simple socle $L$.
 If  $L$  is an alternating group, then a classification of the possible rank $3$ actions appears in \cite{B-1972}, while the rank 3 representations of the classical
groups are listed in \cite{KL-rank3-1982}, and
the rank 3 primitive
groups in which  $L$ is either an exceptional group of Lie type or a sporadic group are listed in \cite{LS-1986}.
(A summary of this classification can also be found in \cite{BM-1994}, which provides also
a list of the smallest possible groups in the various families.) We now inspect the groups case by case to identify the remaining three examples in (1).
Note that the socle $L$ of an almost simple primitive group is not regular, so $L_u\ne1$ for $u\in V(\Gamma)$, and since $\Gamma$ is connected, $L_u^{\Gamma(u)}\ne1$.


\medskip\noindent
\textit{  $L$ exceptional:}
Assume  that  $L$ is an exceptional
simple group of Lie type.  Then $L,L_u$ and the subdegrees $k=|\Gamma(u)|,l=|\Gamma_2(u)|$ are listed in \cite[Table 1]{LS-1986}.
By Lemma \ref{lemma-2atdegree}, we have $k<l$.
In each case, $k$ is not a prime power. Hence $X_u^{\Gamma(u)}$ is a 2-transitive group  on   $\Gamma(u)$ of almost simple type, and
the socle $soc(X_u^{\Gamma(u)})$ is its unique minimal normal subgroup.
Since, as we observed above, $1\ne L_u^{\Gamma(u)}\unlhd X_u^{\Gamma(u)}$,  it follows that
$soc(X_u^{\Gamma(u)})\leq L_u^{\Gamma(u)}$.
Moreover, since $k\neq 28$,
it follows from the classification of the $2$-transitive groups $X_u^{\Gamma(u)}$  that
$soc(X_u^{\Gamma(u)})$ induces a 2-transitive action on $\Gamma(u)$.
Hence  $L_u^{\Gamma(u)}$  acts 2-transitively  on $\Gamma(u)$.
Checking the candidates in \cite[Table 1]{LS-1986}, we see that none of the groups $L_u$ induces a
2-transitive action on degree $k$. Thus $L$ is not exceptional of Lie type.

\medskip\noindent
\textit{ $L$ sporadic:}
For the sporadic
simple groups $L$,
we inspect the groups in \cite[Table 2]{LS-1986}, and find that $(L,L_u)=(M_{22},\mathbb{Z}_2^4.A_6)$ and  $(L,L_u)=(HS,M_{22})$ are  the only two candidates. Moreover,  $(L,L_u)=(M_{22},\mathbb{Z}_2^4.A_6)$,  provides the example  $M_{22}$-graph; and $(L,L_u)=(HS,M_{22})$ yields the Higman-Sims graph.

\medskip\noindent
\textit{ $L$ alternating:}
Suppose  that  $L=A_c$, $c\geq 5$,  as in \cite{B-1972}. Recall that  $n=|\Gamma(u)|$.
Then since $X_u^{\Gamma(u)}$ is 2-transitive on $\Gamma(u)$, it follows from
\cite[Main Theorem]{PraegerWang-1996} that  one of the
following holds:

 \begin{itemize}
\item[(1)]   $X$ is 3-transitive on $V(\Gamma)$;

\item[(2)]   $X\leq S_c$, $n=\frac{c+1}{2}$ and $X_u\cong (S_n\times S_{n-1})\cap X$;

\item[(3)]   $X=S_c$, $c$ is prime, and $X_u=\mathbb{Z}_c:\mathbb{Z}_{c-1}$ is a sharply 2-transitive Frobenius group;

\item[(4)] $X\nleq S_c$, $c=6$, and $X_u$ is the normaliser of a Sylow 5-subgroup;

\item[(5)] $X_u$ is almost simple and primitive in the natural action of $S_c$ of degree $c$, and $X_u^{\Gamma(u)}\cong X_u$.

\end{itemize}

Since $\Gamma$ is not a complete graph, case (1) does not occur.
If  case (2) holds, then by \cite[p.476, I]{B-1972}, $X$ has rank $n$ on $V(\Gamma)$, and so $n=3$ and $c=5$.
In this case, $\Gamma$ is the Petersen graph of girth 5, a contradiction since $\Gamma$ has girth $4$.
Suppose that case (3) holds. Then $|V(\Gamma)|=\frac{|X|}{|X_u|}=(c-2)!$, and the only $2$-transitive action for $X_u$ has degree $c$ (since the valency is not 2), so the valency is $n=c$. Recall that
$|V(\Gamma)|\leq n^2+1$. Hence $(n-2)!\leq n^2+1$, and so $n\leq 6$. Since $n=c$ is a prime power, we get $n=5$,
$|V(\Gamma)|=6$, and $\Gamma\cong \K_6$, which is a contradiction.
In  case (4), we have  $c=6$ and $X_u=N_X(\mathbb{Z}_5)$. Since $X=L.X_u$, it follows that
$|V(\Gamma)|=\frac{|X|}{|X_u|}=|\frac{L.X_u}{X_u}|=|\frac{L}{L\cap X_u}|=\frac{6!}{2.5.2}=36$. Also, the only 2-transitive action of $X_u$ is of degree 5, so the valency is $n=5$, giving
$|V(\Gamma)|>n^2+1$,  a contradiction.
Thus case (5) holds.  Then by
Tables 1-2 of \cite[p.484--485]{B-1972}, we have $c\leq 12$. Since  $X_u$ is  primitive in the natural action of degree $c$, it follows from
\cite[Section 2, I and II]{B-1972} that
there is a unique  graph and it satisfies: $c=9$, $|V(\Gamma)|=120$, $X=A_9$ and $X_u=PSL(2,8):\mathbb{Z}_3\cong {\rm Ree}(3)$. Since
$X_u^{\Gamma(u)}$ is 2-transitive on $\Gamma(u)$ of degree $n$ and $|V(\Gamma)|\leq n^2+1$, it follows that $n=28$. This implies, however, that $X_u$ is transitive on $\Gamma_2(u)$ of degree $120-1-28 = 91$, which is impossible since $|X_u|$ is not divisible by 13.
Thus $L$ is not an alternating group.

\medskip\noindent
\textit{$X$ classical:}
It remains to consider the case where $X$ is a classical almost simple group with a rank 3 action. Such groups are classified in \cite[Theorems 1.1 and 1.2]{KL-rank3-1982}. We consider first the case of \cite[Theorem 1.1]{KL-rank3-1982} which assumes that a group $G$ is semilinear with quasisimple normal subgroup $M= {\rm Sp}(2m-2,q), \Omega^\pm(2m,q), \Omega(2m-1,q)$, or ${\rm SU}(m,q)$, where $m\geq3$ and $q$ is a prime power (and hence $M$ has a possibly non-trivial centre $Z$), so that the socle $L=M/Z$ and the group $X=G/Z\leq \Aut(M/Z)$.
In  cases (i)--(iv) of \cite[Theorem 1.1]{KL-rank3-1982}, the action is on an orbit of totally singular (or isotropic) subspaces, or nonsingular subspaces, and one checks as follows that these do not produce examples: in all cases $X_u$ acts faithfully on each of its orbits $\Gamma(u), \Gamma_2(u)$ and one sees, either group theoretically that the action of $X_u$ on $\Gamma(u)$ is not 2-transitive, or geometrically
that $\Gamma(u)$ contains an edge. (Some additional information from \cite[p.18,Table 6, or the discussion on p.32--36]{BM-1994} is helpful.)
The remaining quasisimple groups $M$ in cases (v)--(x) of \cite[Theorem 1.1]{KL-rank3-1982}
are listed in Table~\ref{t:diam2}, together with the stabiliser $M_u$, and the subdegrees (lengths of the $G_u$-orbits in $V(\Gamma)$) many of which were computed using MAGMA \cite{Magma-1997}. In most cases we determine that $G_u$ does not act 2-transitively on either $\Gamma(u)$ or
$\Gamma_2(u)$  (often using MAGMA to confirm) - we denote this fact by an entry `not $2$-trans.'  in the column headed `Comments' of Table~\ref{t:diam2}. In one case, case (vi), we find the 2-arc-transitive Hoffman-Singleton graph, but it is excluded here since it has girth 5, while we are assuming that $\Gamma$ has girth $4$.

\begin{table}
\caption{Cases (v)--(x) of \cite[Theorem 1.1]{KL-rank3-1982} for the Proof of  Theorem \ref{diam2-theo-1}}
\begin{tabular}{llll}
\hline
$M$ &Stabiliser $M_u$& Subdegrees & Comments  \\
\hline\hline 
$SU(3,3)$    & $PSL(3,2)$   & $1, 14, 21$ 	&  not $2$-trans. \\
$SU(3,5)$    & $3.A_7$       & $1, 7, 42$ 	&  Hoffman-Singleton not girth $4$\\
$SU(4,3)$    & $4.PSL(3,4)$& $1, 56, 105$ &  not $2$-trans.  \\
$Sp(6,2)$    & $G_2(2)\cong U_3(3):2$		&$1, 56, 63$ &  not $2$-trans.  \\
$\Omega(7,3)$    & $G_2(3)$     & $1, 351, 728$ 		&  not $2$-trans. \\
$SU(6,2)$    & $3.PSU(4,3):2$& $1, 56, 105$ &  not $2$-trans.  \\
\hline
\end{tabular}
\label{t:diam2}
\end{table}

Finally, we treat the groups from   \cite[Theorem 1.2]{KL-rank3-1982}. Here
$L=PSL(n,q)\leq X\leq \Aut(L)$. In case (i) of  \cite[Theorem 1.2]{KL-rank3-1982}, $n\geq4$
and the $X$-action on $V(\Gamma)$ is equivalent to its action on lines of the projective space $PG(n-1,q)$. For a line $u$, the sets $\Gamma(u), \Gamma_2(u)$ are the sets of lines that either
intersect $u$ in a single point, or are disjoint from $u$. In either case it is easy to see that the corresponding graph contains a triangle, and hence we obtain no girth 4 example. Of the remaining groups in cases (ii)--(iv) of  \cite[Theorem 1.2]{KL-rank3-1982}, we do not need to consider
$L=PSL(2,4), PSL(2,9)$ or $PSL(4,2)$ as these groups are isomorphic to alternating groups, which were treated above. With these exclusions,
the remaining groups $L$ in these cases are listed in
 Table~\ref{t:diam2-2}, together with the stabiliser $L_u$, and the subdegrees. In all but one
 case $X_u$ does not act 2-transitively on either $\Gamma(u)$ or
$\Gamma_2(u)$ - denoted as above by `not $2$-trans.'  in the column headed `Comments'. In the exceptional case, case (iii) of  \cite[Theorem 1.2]{KL-rank3-1982} with $L=PSL(3,4)$, the stabiliser is 2-transitive on its orbit of length 10, and we obtain the 2-arc-transitive Gewirtz graph of girth $4$, as in Part (1).  This completes the proof  of  Theorem \ref{diam2-theo-1}. \qed

\begin{table}
\caption{Cases (ii)--(iv) of \cite[Theorem 1.2]{KL-rank3-1982} for the Proof of  Theorem \ref{diam2-theo-1}}
\begin{tabular}{llll}
\hline
$L$ &Stabiliser $L_u$& Subdegrees & Comments  \\
\hline\hline 
$PSL(2,8)$    & $D_{14}$   & $1, 14, 21$ 	&  not $2$-trans. \\
$PSL(3,4)$    & $A_6$       & $1, 10, 45$ 	&  Gewirtz graph \\
$PSL(4,3)$    & $PSp(4,3):2$& $1, 36, 80$ &  not $2$-trans.  \\
\hline
\end{tabular}
\label{t:diam2-2}
\end{table}

\bigskip

As a corollary of Theorem \ref{diam2-theo-1}, we observed a group theoretic characterisation of the Petersen graph.
A primitive group on a set $X$ is \emph{$2$-primitive} if for every $u\in X$, the stabiliser of $u$  is primitive on $X\setminus \{u\}$. In particular, each 3-transitive group
is 2-primitive.

\begin{corol}\label{2primitive-faith-1}
Let $\Gamma$ be a  connected $(G,2)$-arc-transitive graph of girth   $4$ or $5$,
and let  $u\in V(\Gamma)$. Assume that the $G_{u}$-action  on $\Gamma(u)$ is $2$-primitive and unfaithful. Then
$\Gamma$ is either the Petersen graph or a complete bipartite graph.

\end{corol}

\proof Let  $(u,v)$ be an arc of  $\Gamma$.
Since $G_{u}$ is   $2$-primitive on $\Gamma(u)$,
$G_{uv}$ acts primitively on $\Gamma(u)\setminus \{v\}$.
Let $K$ be the kernel of $G_v$ acting on $\Gamma(v)$. Then    $K\unlhd G_v$ and $K \unlhd G_{uv}$.
Since $G_{uv}$ is primitive on $\Gamma(u)\setminus \{v\}$,   either $K$ fixes all the vertices of
$\Gamma(u)\setminus \{v\}$ or $K$ is transitive on
$\Gamma(u)\setminus \{v\}$. If $K$ fixes all the vertices in
$\Gamma(u)\setminus \{v\}$, then   $K=1$ as $\Gamma$ is connected, and so $G_{u}$ is faithful on $\Gamma(u)$, contradicting our assumption. Hence $K$ is transitive on
$\Gamma(u)\setminus \{v\}$.  Note that $K\leq G_{uvw}$ for each   2-geodesic $(u,v,w)$. Thus  $G_{uvw}$
is transitive on $\Gamma(u)\setminus \{v\}$.
Suppose first that $\Gamma$ has girth 4. Then since $G$ is transitive on $2$-arcs, $\Gamma$ has a cycle of length 4
containing $u, v, w$, and so $w$ is adjacent to at least one vertex of $\Gamma(u)\setminus \{v\}$.
Since $G_{uvw}$
is transitive on $\Gamma(u)\setminus \{v\}$, it follows that
$\Gamma(u)=\Gamma(w)$, and since this holds for all $w$, it follows that    $\Gamma$ is a complete bipartite graph.
Suppose then that $\Gamma$ has girth 5. Then, for each of the vertices
$x\in\Gamma(u)\setminus\{v\}$ there exists a vertex $y(x)\in \Gamma(w)\setminus\{v\}$ such that $(u,v,w,y(x),x)$ is a $5$-cycle. Since $\Gamma$ contains no 4-cycles, distinct $x, x'$ correspond to distinct $y(x), y(x')$ and it follows that $\Gamma(w)\setminus\{v\}\subseteq \Gamma_2(u)$
so $\Gamma$ has diameter $2$.
Then by Theorem \ref{diam2-theo-1}, $\Gamma$  is one of the following three graphs: $C_5$,
the Petersen graph,  and the Hoffman-Singleton graph. However, if $\Gamma$ is  $C_5$ or the Hoffman-Singleton graph, then the vertex stabilisers (in the full
automorphism group) act faithfully on their neighbours, a contradiction. Hence
$\Gamma$  is the Petersen graph.
\qed

\section{Proof of Theorem \ref{girth45-redtheo-1}}

In this section, we investigate the class of 3-geodesic-transitive  graphs of girth 4 or 5 which are normal covers of a 2-arc-transitive graph of diameter at most 2.
The first subsection determines such covers of complete graphs. We use the parameters $a_i, b_i, c_i$ introduced in Definition~\ref{intersectionarray}.

\subsection{Covers of complete graphs}

\begin{lemma}\label{quot-comp-2}
Let $\Gamma$ be a connected  $(G,3)$-geodesic-transitive  graph of girth $4$ or $5$.
Suppose that    $\Gamma$ is a $G$-normal cover of the  complete graph $\K_r$ for some
$r\geq 3$, relative to a normal subgroup $N$ of $G$.
Then $\Gamma$ is a distance-transitive antipodal cover of $\K_r$ of diameter $3$,
and either
\begin{enumerate}
\item[(i)]$\Gamma=\K_{r,r}- r\K_2$ for some $r\geq 4$, or
\item[(ii)] $\Gamma$ is the distance $2$ graph $[HoS]_2$ of the Hoffman--Singleton graph, with $r=7$.
\end{enumerate}
\end{lemma}

\proof
Let $\mathcal{B}=\{B_1,\ldots,B_r\}$ be the set of $N$-orbits in $V(\Gamma)$, so that
$\Sigma := \Gamma_{\mathcal{B}} \cong \K_r$ with $r\geq 3$.
Since $\Gamma$ is a cover of $\Sigma$, the valency of $\Gamma$ is $r-1\geq 2$.
Since $\Gamma$ has girth 4 or 5, $\Gamma$ is not a complete graph, and in particular $|B_i|\geq2$. For the same reason,
if $r=3$, then $\Gamma$ has valency 2 and so $\Gamma$ is $C_4$ or $C_5$, neither of which is a cover of $\K_3$. Hence $r\geq 4$.
Let $(u_1,u_2,u_3)$ be a 2-geodesic of $\Gamma$ with $u_i\in B_i$, for $i=1, 2$. Since  $\Gamma$ is a cover of $\Sigma$,  $u_3\notin B_1$, so we may suppose that $u_3\in B_3$. As  $\Sigma$ is a complete graph, it follows that
$B_1$ and $B_3$ are adjacent in $\Sigma$, and so $u_3$ is adjacent in $\Gamma$ to a vertex $u_1'\in B_1$. Now $u_1'\ne u_1$ since $\Gamma$ has girth at least 4, and since $\Gamma$
is   $G$-arc-transitive, $B_1$ does not contain an edge of $\Gamma$ and hence $u_1'\not\in\Gamma(u_1)$.  Thus
$u_1'\in \Gamma_2(u_1)\cup \Gamma_3(u_1)$.
If $u_1'\in \Gamma_2(u_1)$, then there exists a vertex $v$ such that $(u_1,v,u_1')$ is a 2-geodesic. This implies that  $\{u_1,u_1'\}\subseteq \Gamma(v)\cap B_1$, so  $|\Gamma(v)\cap B_1|\geq 2$,
contradicting the fact that $\Gamma$ is a cover of $\Sigma$. Thus
$u_1'\in  \Gamma_3(u_1)$, so $u_1'\in  \Gamma_3(u_1)\cap B_1$. Since  $\Gamma$ is  $(G,3)$-geodesic-transitive,   the stabiliser $G_{u_1}$ is transitive on $\Gamma_3(u_1)$ and fixes  $B_1$ setwise, so $\Gamma_3(u_1)\subset B_1$.
Furthermore, as  $B_1$ does not contain an edge of $\Gamma$,  the induced subgraph $[\Gamma_3(u_1)]$ is an empty graph.

Since the diameter of $\Gamma$ is  at least 3, we have  $b_2=|\Gamma(u_3)\cap \Gamma_3(u_1)|\geq 1$.
On the other hand, as  $\Gamma$  is a cover of  $\Sigma=\K_r$ and $\Gamma_3(u_1)\subset B_1$,  it follows that $|\Gamma(u_3)\cap B_1|=1$. Thus  $b_2=1$.
Since $\Gamma$ is   $(G,3)$-geodesic-transitive, it is easy to see  that $b_3\leq b_2=1$, and so
\begin{center}
$b_3=0$ or $1.$
\end{center}

Suppose  that $b_3=1$. Then  $\Gamma$ has diameter at least 4. As    $[\Gamma_3(u_1)]$ is an empty graph, we obtain $a_3=0$, and so  $c_3=(r-1)-b_3=r-2$. Let $(u_1,u_2,u_3,u_4,u_5)$ be a 4-geodesic of $\Gamma$. Then $u_4\in B_1$ and $u_5\notin B_1$. Hence $|\Gamma(u_5)\cap B_1|=1$, and so $\Gamma(u_5)\cap B_1=\{u_4\}$,  and $|\Gamma(u_5)\cap \Gamma_3(u_1)|=1$ (since $\Gamma_3(u_1)\subset B_1$). However,
since  $(u_2,u_3,u_4,u_5)$ is  a 3-geodesic and  $c_3=r-2$,  we have $| \Gamma(u_5)\cap \Gamma_2(u_2)|=r-2$. We note that  $ \Gamma(u_5)\cap \Gamma_2(u_2)\subseteq \Gamma(u_5)\cap \Gamma_3(u_1)$ (since $u_5\in\Gamma_4(u_1)$). Thus  $r-2 =| \Gamma(u_5)\cap \Gamma_2(u_2)|\leq 1$,  contradicting the fact that $r\geq 4$.

Thus  $b_3=0$, and  $\Gamma$ has diameter 3. Hence $\Gamma$  is   $G$-geodesic-transitive.
Since $\Gamma$ has valency $r-1$ and girth  $4$ or 5, it follows that $b_0=r-1,b_1=r-2$ and $c_1=1$, and we have shown that $b_2=1$. Furthermore, $2\leq c_2\leq r-2$ if the girth is $4$  and $c_2=1$ if the girth is $5$.
Since  $\Gamma$ has diameter 3 and  $[\Gamma_3(u_1)]$ is an empty graph (since $\Gamma_3(u_1)\subset B_1$), it follows that  $a_3=0=b_3$, and so $c_3=r-1$. Hence
$\Gamma$ has intersection array
$$
(r-1,r-2,1;1,c_2,r-1).
$$

Further, for distinct $z,z'\in \Gamma_3(u_1)$, the distance
$d_{\Gamma}(z,z')\neq 1$ since $a_3=0$ and
$d_{\Gamma}(z,z')\neq 2$ since $b_2=1$, and hence $d_{\Gamma}(z,z')=3$. This implies that
$\Gamma$ is a $G$-distance-transitive  antipodal cover of $\Sigma$, with antipodal blocks of size $|B_i|=1+|\Gamma_3(u_1)|=1+(r-2)/c_2$.
Thus, $\Gamma$ is listed in  \cite[Main Theorem]{GLP}.
Moreover, since $\Gamma$ is  $(G,3)$-geodesic-transitive  with girth 4 or 5, it follows that for each vertex $u$, $G_u$ is 2-transitive on $\Gamma(u)$, and hence $G$ induces a 3-transitive group on $V(\Sigma)$.
We inspect the candidates in \cite[Main Theorem]{GLP}.

The  graph  in case  (1) of \cite[Main Theorem]{GLP} is $\Gamma= \K_{r,r}-r\K_2$,
as in part (i), and the graph in  case  (2) is  $\Gamma = [HoS]_2$ with $r=7$  (see also
\cite[p.223]{BCN}). Both are geodesic-transitive.  None of the groups $G$ in part (3) (a)--(e)
 of \cite[Main Theorem]{GLP} induce a 3-transitive group on $V(\Sigma)$.

Let $\Gamma$ be a graph in case (4)(a). Then $\Gamma$ is also in  \cite[Propositon 12.5.3]{BCN}, its intersection array is
$(q,q-\frac{q-1}{r}-1,1;1,\frac{q-1}{r},q)$, so
$\Gamma$ has girth 3, a contradiction.
For graphs in  Main Theorem (4)(b) and (5) of \cite{GLP}, $G/N \leq \Aut(\Gamma)$ does not induce a
3-transitive action on $V(\Sigma)$, where $N$ is the kernel of $G$ acting on $V(\Sigma)$.

Finally, let $\Gamma$ be a graph in case   (6). Then $\Gamma$ is described in Example 3.6 and  Section 6 of
\cite{GLP}. By Lemma 6.1 of \cite{GLP}, $\Gamma$ is a normal Cayley graph of some $p$-group $P$ and $\Gamma$ is as in Proposition 6.2 or 6.3 of \cite{GLP}. Let $u=1_P$. Since $G_u$ is 2-transitive on $\Gamma(u)$, we have  $p=2$, and
hence the graphs in  \cite[Proposition 6.3]{GLP} do not occur.
Thus $\Gamma$ is as in \cite[Proposition 6.2]{GLP}. Here $r=2^{2b}$ and $|B_i|=2^a$ for some positive integers $a, b$,
and $G$ induces an affine $3$-transitive group on $\Sigma=K_r$
contained in $[2^{2b}].Sp(2b,2)$. Such a group is $3$-transitive only when $b=1$, and in this case $2^a=|B_i|=1+(r-2)/c_2=1+2/c_2$,
so that $c_2=2$ and $|B_i|=2$. It follows that  $\Gamma\cong \H(3,2)\cong Q_3\cong \K_{4,4}-4\K_2$, as in (i).
 \qed

\subsection{Covers of the Petersen graph and $\Box_5$, and a general lemma}

First we  determine the   $3$-geodesic-transitive covers of
the Petersen graph and the folded 5-cube.

\begin{lemma}\label{cover-petersen}
Let $\Gamma$ be a connected  $3$-geodesic-transitive  graph of girth $4$ or $5$. Then the following hold.

{\rm (1)} $\Gamma$ is not a cover of $ C_r$ for $r=4,5$.

{\rm (2)} If   $\Gamma$  is a cover of the Petersen graph, then $\Gamma$ is the dodecahedron.

{\rm (3)} If   $\Gamma$ is a cover of $\Box_5$, then $\Gamma$ is either $\H(5,2)$ or the Armanios-Wells graph.

\end{lemma}
\proof  (1) If  $\Gamma$ is a cover of $ C_r$ where  $r=4,5$, then $\Gamma$ has valency 2 and girth 4 or 5, and so $\Gamma\cong C_4$ or $C_5$ which have  no 3-geodesics, a contradiction.

(2) Suppose that  $\Gamma$  is a cover of the Petersen graph.
Then $\Gamma$ has valency 3 and girth 5, and so  $c_2=1$ and $a_2=1$ or 2.
Since $\Gamma$ has diameter at least $3$, we must have  $a_2=1$ and  $b_2=1$. Then by Remark \ref{3gt-dt-1}, $\Gamma$ is a geodesic-transitive graph of valency 3, so $\Gamma$ is listed in  \cite[p.221]{BCN}. By inspecting the candidates, we conclude that $\Gamma$ is the dodecahedron.

(3) Suppose that  $\Gamma$ is a cover of $\Box_5$.  Since $\Box_5$ has valency 5 and girth 4, it follows that $\Gamma$ has valency 5 and  girth 4 or 5.
Then by Theorem 1.1 of \cite{JW-2018},  $\Gamma$ is either $\H(5,2)$ or
the Armanios-Wells graph. \qed

The following lemma will be used frequently in the rest of our analysis.

\begin{lemma}\label{cover-group-1}
Let $\Gamma$ be a connected  $(G,3)$-geodesic-transitive graph of  girth $4$ or $5$.
Suppose that $G$ has an intransitive  normal subgroup $N$ with at least $3$ orbits
on $V(\Gamma)$ such that
$\Gamma_N$ has   diameter $2$.  Let  $(B_0,B_1,B_2)$ be a
$2$-geodesic of $\Gamma_N$. Then there exist $u_i\in B_i$, for each $i$
such that $(u_0,u_1,u_2)$ is a  $2$-geodesic of $\Gamma$, and for each such $2$-geodesic
the following statements hold.

\begin{itemize}
\item[(1)] $\Gamma_N$ has girth at least $4$.

\item[(2)] $G_{B_0}=NG_{u_0}$, $G_{B_0B_1}=NG_{u_0u_1}$ and $G_{B_0B_1B_2}=NG_{u_0u_1u_2}$.

\item[(3)] If   $B\in \Gamma_N(B_0)\cap \Gamma_N(B_2)$ such that $B\neq B_1$, then
$\Gamma(u_2)\cap B \subseteq  \Gamma(u_0)\cup \Gamma_3(u_0)$. Moreover if
$\Gamma(u_2)\cap B \cap  \Gamma_3(u_0)\ne\emptyset$, then  all vertices  of
$\Gamma_3(u_0)$ lie in blocks of $\Gamma_N(B_0)$ and $|\Gamma_3(u_0)|$ is divisible by the valency
$|\Gamma(u_0)|$.

\item[(4)] If $|\Gamma(u_0)\cap \Gamma(u_2)|=|\Gamma_N(B_0)\cap\Gamma_N(B_2)|$ (that is if
$\Gamma$ and $\Gamma_N$ has `the same $c_2$'), then $G_{B_0B_2}=NG_{u_0B_2}= NG_{u_0u_2}$.

\end{itemize}

\end{lemma}
\proof (1)  Since $\Gamma$ has girth 4 or 5, $\Gamma\not\cong \K_{m[b]}$ for any $m\geq3, b\geq2$.
Hence by  Lemma~\ref{2gt-covergt-1}, $\Gamma$ is a cover of $\Gamma_N$,  $N$ is semiregular on $V(\Gamma)$ (so each $N$-orbit has $|N|$ vertices), and $\Gamma_N$ is $(G/N, 2)$-geodesic-transitive.
Since $\Gamma_N$ has diameter 2 it follows from Lemma~\ref{girth45-lem-1} that $\Gamma_N$ has girth at least 4.

(2) Also, since $\Gamma$ covers $\Gamma_N$, each $u_0\in B_0$ determines a unique $u_1\in B_1$ and
$u_2\in B_2$ such that  $(u_0,u_1,u_2)$ is a $2$-arc of $\Gamma$, and indeed it must be
a $2$-geodesic of $\Gamma$ since  $\Gamma_N$ has girth at least 4.
Now $B_0\in V(\Gamma_N)$ is the $N$-orbit containing $u_0$, that is,
$B_0 = u_0^N$, and $B_0$ is a block of imprimitivity for $G$ in $V(\Gamma)$.
Thus both $G_{u_0}\leq G_{B_0}$ and $N<G_{B_0}$, and as $N$ is transitive on $B_0$, it follows that   $G_{B_0}=NG_{u_0}$.

For $i=1, 2$, consider the transitive action of $G_{B_0\dots B_i}$ on $B_0$. The stabiliser of $u_0$ is $H:=G_{B_0\dots B_i}\cap G_{u_0}$. Since $u_1$ is the unique point of $B_1$ adjacent to $u_0$, $H$ also fixes $u_1$; and similarly, if $i=2$,  since $u_2$ is the unique point of $B_2$ adjacent to $u_1$, $H$ also fixes $u_2$. Thus $H\leq G_{u_0\dots u_i}$, and conversely $G_{u_0\dots u_i}\leq G_{B_0\dots B_i}\cap G_{u_0}$. Hence the stabiliser $H= G_{u_0\dots u_i}$.  Since the subgroup $N$ of
$G_{B_0\dots B_i}$ is transitive on $B_0$, it follows that $G_{B_0\dots B_i}=NG_{u_0\dots u_i}$.

(3) Let $B\in \Gamma_N(B_0)\cap \Gamma_N(B_2)$ such that $B\neq B_1$ and let $u\in \Gamma(u_2)\cap B$.
Since $(u_0,u_1,u_2)$ is a   $2$-geodesic of $\Gamma$, the distance $d_\Gamma(u_0,u_2)=2$ and so $d_\Gamma(u_0,u)\leq 3$. If $d_\Gamma(u_0,u)=2$, then some element of $G_{u_0}$ maps $u$ to $u_2$, and hence maps $B$ to $B_2$, which is impossible since $B, B_2$ have distances 1, 2 from $B_0$ in $\Gamma_N$, respectively.
Thus  $u\in \Gamma(u_0)\cup \Gamma_3(u_0)$. Finally suppose that $d_\Gamma(u_0,u)= 3$. Since  $\Gamma$ is  $(G,3)$-geodesic-transitive, $G_{u_0}$ is transitive on $\Gamma_3(u_0)$, and so all points of
$\Gamma_3(u_0)$ lie in blocks of $\Gamma_N(B_0)$. Since $G_{u_0}$ is transitive on $\Gamma_N(B_0)$,
it follows that $|\Gamma_N(B_0)|=|\Gamma(u_0)|$ divides $|\Gamma_3(u_0)|$.

(4)  Suppose that $\Gamma$ and $\Gamma_N$ have the same value of $c_2$. Since $G$ is transitive on the 2-geodesics of both $\Gamma$ and $\Gamma_N$, we have $c_2= |G_{B_0B_2}:G_{B_0B_1B_2}|=|G_{u_0u_2}:G_{u_0u_1u_2}|$. Hence using this and part (2), we have
\[
\frac{|G_{B_0B_2}|}{|G_{u_0u_2}|} =
\frac{|G_{B_0B_1B_2}|}{|G_{u_0u_1u_2}|} = |N|,
\]
and it follows that $G_{B_0B_2}=NG_{u_0u_2}$. Since $NG_{u_0u_2}\subseteq
N G_{u_0B_2}\subseteq G_{B_0B_2}$, equality holds and part (4) is proved.
\qed

\subsection{Covers of the Higman-Sims graph}
\medskip
Recall that the Higman-Sims graph  is strongly regular with parameters $(100,22,0,6)$.
The following lemma determines the unique $3$-geodesic-transitive cover of the Higman-Sims graph.

\begin{lemma}\label{cover-hs-1}
Let $\Gamma$ be a connected  $(G,3)$-geodesic-transitive graph of girth $4$ or $5$ such that $\Gamma$ is a
$G$-normal cover of the Higman-Sims graph $\Sigma$. Then  $\Gamma$ is the standard double cover of $\Sigma$, and
$G=H\times \mathbb{Z}_2$, where $H$ is the Higman-Sims group $HS$ or $HS.\mathbb{Z}_2$.
\end{lemma}

\proof By assumption there is a non-trivial normal subgroup $N\unlhd G$ such that $\Gamma$ is a cover of  the Higman-Sims graph $\Sigma=\Gamma_N$.  Since $\Sigma$ is strongly regular with parameters $(100,22,0,6)$, the graph
$\Gamma$ has valency 22 and $c_2\leq 6$.
Let $u\in V(\Gamma)$, and let $B = u^N$ be the $N$-orbit (vertex of $\Sigma$) containing $u$, so   $G_B=NG_u$ (Lemma~\ref{cover-group-1} (2)). By Lemma~\ref{2gt-covergt-1},
$N$ is semiregular on $V(\Gamma)$,  the subgroup induced by $G$ on $\Sigma$ is $G/N$ and this group acts transitively on the $2$-geodesics of $\Sigma$. In particular, $G_B$ acts 2-transitively on $\Sigma(B)$, and it follows that  $G/N\cong HS$ or $HS.\mathbb{Z}_2$, and $G_u\cong G_B/N\cong M_{22}$ or $M_{22}.\mathbb{Z}_2$, respectively (see \cite[p.39, 80]{Atlas}, for example).

Let  $(B_0,B_1,B_2)$ be a 2-geodesic of $\Sigma$. Since $\Gamma$ covers $\Sigma$, there are vertices $u_i\in B_i$, for $i=0, 1, 2$
such that $(u_0,u_1,u_2)$ is a  $2$-geodesic of $\Gamma$. It follows from the parameters of $\Sigma$ that the
set $\Delta:= \Sigma(B_0)\cap \Sigma(B_2)$
has size 6 and contains $B_1$. Moreover, exactly $c_2\geq 1$ of the $\Sigma$-vertices  in $\Delta$ (including $B_1$) contain a vertex of $\Gamma(u_0)\cap \Gamma(u_2)$.  In particular, this $c_2$-subset of   $\Delta$ is fixed setwise by $G_{u_0u_1u_2}$. By Lemma~\ref{cover-group-1} (2), $G_{B_0B_1B_2}=NG_{u_0u_1u_2}$, and hence $G_{B_0B_1B_2}$ leaves invariant a $c_2$-subset of $\Delta$ containing $B_1$.

Considering the $G$-action on $\Sigma$, we have (see \cite[p.80]{Atlas}) $G_{B_0B_2}/N = \mathbb{Z}_2^4.A_6$ or $\mathbb{Z}_2^4.S_6$, and so $G_{B_0B_2}$ induces $A_6$ or $S_6$ on $\Delta$, so that $G_{B_0B_1B_2}$ is transitive on the $5$-subset $\Delta\setminus\{B_1\}$. We conclude from the observation in the previous paragraph that $c_2=1$ or $c_2=6$.  If $c_2=6$, then  it follows from work of Cameron, see \cite[Part (II) on p.4]{CP-1983}, that $\Gamma$ is the standard double cover of the Higman-Sims graph as in Definition \ref{double-def}, and the assertions of the lemma hold.
Assume to the contrary that $c_2=1$, so $\Gamma$ has girth $5$. We shall obtain a contradiction.

Let $B_3\in (\Sigma(B_0)\cap \Sigma(B_2))\setminus\{B_1\}$
and $B_4\in \Sigma_2(B_0)\cap \Sigma(B_2)$.
Let $u_i\in B_i$ ($i=3,4$) such that $(u_0,u_1,u_2,u_3)$
and $(u_0,u_1,u_2,u_4)$ are 3-arcs. Then as $c_2=1$,   the two vertices $u_3$ and $u_4$ are in $\Gamma_2(u_0)\cup \Gamma_3(u_0)$.
By Lemma~\ref{cover-group-1}(3), it follows that $u_3\in\Gamma_3(u_0)$, and since this holds for each of the $5$ blocks
$B_3\in  (\Sigma(B_0)\cap \Sigma(B_2))\setminus\{B_1\}$, it follows that $b_2 = |\Gamma_3(u_0)\cap\Gamma(u_2)|\geq 5$.
Moreover, by  Lemma~\ref{cover-group-1}(3),  all points of
$\Gamma_3(u_0)$ lie in blocks of $\Sigma(B_0)$, and therefore $u_4\in\Gamma_2(u_0)$.
Since this holds for each of the $16$ blocks $B_4\in \Sigma_2(B_0)\cap \Sigma(B_2)$, it follows that
$a_2 = |\Gamma_2(u_0)\cap\Gamma(u_2)|\geq 16$ for $\Gamma$. Thus $22=|\Gamma(u_2)|=a_2+b_2+c_2\geq 16+5+1 = 22$, and so  $a_2=16$ and $b_2=5$.

Since $B_3\in \Sigma(B_0)$, there exists a vertex  $u_0'\in B_0$ such that
$(u_3,u_0')$ is  an arc of $\Gamma$. Further $u_0'\neq u_0$, as $\Gamma$ has
girth 5. The sequence $(u_0,u_1,u_2,u_3, u_0')$ is a 4-arc, so $1\leq
d_\Gamma(u_0,u_0')\leq 4$. Further $G_{u_0}$ fixes $\Sigma(B_0)$ and $\Sigma_2(B_0)$ setwise, and there are blocks
in $\Sigma(B_0)\cup \Sigma_2(B_0)$ that contain vertices from $\Gamma_i(u_0)$ for $i=1, 2$ or $3$ (such as $u_1, u_2, u_3$ respectively).  Since
$\Gamma$ is $(G,3)$-geodesic-transitive,  $G_{u_0}$ is transitive on $\Gamma_i(u_0)$ for each $i=1, 2, 3$, and hence
$d_\Gamma(u_0,u_0')=4$, and $\Gamma$ has diameter at least $4$..

Consider $\Phi:=\Gamma_2(u_1)\cap \Gamma(u_{3})$. Then  $\Phi\subseteq \Gamma_2(u_0)\cup \Gamma_3(u_{0})$. Since $d_\Gamma(u_1,u_3)=2$, it follows that $|\Phi|=a_2=16$.
Let $u\in \Phi$ and let $B$ be the block containing  $u$. If  $u\in \Gamma_3(u_{0})$,
then $B$ would lie in $\Sigma(B_0)$ and $(B_0, B_3, B)$ would be a triangle in $\Sigma$,
contradicting the fact that $\Sigma$ has girth $4$. Thus $\Phi\cap\Gamma_3(u_0)=\emptyset$, and
hence    $\Phi\subseteq \Gamma_2(u_0)$, so $c_3=|\Gamma_2(u_0)\cap \Gamma(u_3)|\geq |\Phi|=16$. Also $c_3\leq 21$ since $\Gamma$ has diameter at least $4$.
Counting the edges between $\Gamma_2(u_0)$ and $\Gamma_3(u_0)$, we find $|\Gamma_3(u_0)|\cdot c_3 = |\Gamma_2(u_0)|\cdot b_2
= 22\cdot 21\cdot 5$, and as $16\leq c_3\leq 21$ we conclude that $c_3=21$.
Since $\Gamma$ has valency 22 and diameter at least 4, this implies that $b_3= 1$, $a_3=0$.
 Thus by Remark \ref{3gt-dt-1} (2), $\Gamma$ is $G$-distance-transitive.
Furthermore, $|\Gamma_3(u_0)|=22\cdot 5$, and so the number of edges between $\Gamma_3(u_0)$ and $\Gamma_4(u_0)$
is $22\cdot 5 = |\Gamma_3(u_0)|\cdot b_3 = |\Gamma_4(u_0)|\cdot c_4$. Thus $c_4$ divides $22\cdot 5$ and by \cite[Proposition 20.4]{Biggs-1},
$c_4\geq c_3=21$, so $c_4=22,  |\Gamma_4(u_0)|=5$, and $\Gamma$ has diameter 4. The set of $N$-orbits is a block system for $G$ which is not a bipartition, and it follows from \cite[Theorem 2]{smith} that it is a set of antipodal blocks and $\Gamma$ is antipodal.
Thus $B_0=\{u_0\}\cup\Gamma_4(u_0)$ and $G_{B_0}^{B_0}$ is 2-transitive on $B_0$ of degree $6$. However as $N$ is semiregular on $V(\Gamma)$, $N^{B_0}$ is a regular normal subgroup of $G_{B_0}^{B_0}$, which is a contradiction since $|B_0|=6$ is not a prime power.
\qed

\subsection{Covers of the Gewirtz graph}

\medskip
Next we consider covers of the Gewirtz graph. Our first step identifies two substantial cases to be analysed.

\begin{lemma}\label{cover-g-1A}
Let $\Gamma$ be a connected  $(G,3)$-geodesic-transitive graph of girth $4$ or $5$.
Suppose that $\Gamma$ is a $G$-normal cover of the Gewirtz graph $\Sigma$, that is, $G$ has a normal subgroup $N$ such that $\Gamma_N\cong \Sigma$. Then  the following hold:

 \begin{itemize}
\item[(1)] $\Gamma$ has girth $4$ and $c_2=2$;

\item[(2)] $PSL(3,4)\leq G/N\leq PSL(3,4).\mathbb{Z}_2^2$ and, for $u_0\in V(\Gamma)$, either

\begin{enumerate}
\item[(a)]  $G_{u_0}=PSL(2,9)$ or $P\Sigma L(2,9)$, and $(a_2, b_2)=(4,4)$; or

\item[(b)]  $G_{u_0}=PGL(2,9), M_{10}$ or $P\Gamma L(2,9)$, and $(a_2, b_2)=(0,8)$.
\end{enumerate}
\end{itemize}

\end{lemma}

\proof We identify $\Gamma_N=\Sigma$. Since $\Sigma$ is strongly regular
with parameters $(56,10,0,2)$, it follows that $\Gamma$ has valency 10 and $c_2\leq 2$. So either $\Gamma$ has girth $4$ and $c_2=2$, or $\Gamma$ has girth $5$ and $c_2=1$.
Let  $(B_0,B_1,B_2)$ be a 2-geodesic of $\Sigma$ and let $B_3\in \Sigma(B_0)\cap \Sigma(B_2)$
such that $B_3\neq B_1$. Let $u_i\in B_i$, for each $i$, such that $(u_0,u_1,u_2,u_3)$ is a 3-arc of $\Gamma$.
Since $\Gamma$ has girth $4$ or $5$, it follows that $(u_0,u_1,u_2)$ is a 2-geodesic.
Let  $B_4\in \Sigma_2(B_0)\cap \Sigma(B_2)$ and let $u_4\in B_4$  so that $(u_0,u_1,u_2,u_4)$ is  a 3-arc.

Suppose for a contradiction that  the girth of $\Gamma$ is 5, that is,  $c_2=1$. By Lemma~\ref{cover-group-1}(3), the  vertex  $u_3\in \Gamma_3(u_0)$ so $(u_0,u_1,u_2,u_3)$ is a 3-geodesic,
and $\Gamma_3(u_0)$ is contained in the union of the blocks of $\Sigma(B_0)$. This implies that $u_4\notin \Gamma_3(u_0)$, and so  $u_4\in \Gamma_2(u_0)$. Since this holds for each choice of $B_4$ and  $|\Sigma_2(B_0)\cap \Sigma(B_2)|=8$,
it follows that $a_2=|\Gamma_2(u_0)\cap\Gamma(u_2)| \geq 8$, and since $\Gamma$ has valency 10
and $\Gamma(u_2)$ contains also $u_1, u_3$, we must have $a_2=8$, $b_2=1$ and $c_2=1$.
Thus  by Remark \ref{3gt-dt-1}(1), $\Gamma$ is a $G$-distance-transitive graph of valency 10. However, by inspecting the candidates in \cite[p.224]{BCN}, such a graph does not exist.
Hence $\Gamma$ has girth 4 and so $c_2=2$, and Part (1) holds.
Also, since  $\Sigma(B_0)\cap \Sigma(B_2)=\{B_1, B_3\}$, we have  $\Gamma(u_0)\cap \Gamma(u_2)=\{u_1, u_3\}$,
$(u_0,u_1,u_2,u_3)$ is a 4-cycle, and  $(u_0,u_1,u_2,u_4)$ is  a 3-geodesic  for at least one choice of $B_4, u_4$.

 The group $G$ induces a subgroup $G/N$ of $\Aut(\Sigma) = PSL(3,4).\mathbb{Z}_2^2$, and by \cite[p.23]{Atlas}, the only vertex-transitive proper subgroups of $\Aut(\Sigma)$ contain $\PSL(3,4)$. Thus $PSL(3,4)\leq G/N$. By Lemma~\ref{cover-group-1} (2), $G_{u_0}\cong G_{B_0}/N$ and so
$PSL(2,9)\leq G_{u_0}\leq PSL(2,9).\mathbb{Z}_2^2$ (see \cite[p.23]{Atlas}), as in Part (2).
Since $\Gamma$ is $(G,3)$-geodesic-transitive,  $G_{u_0u_1u_2}$ is transitive on the $b_2$-subset
$\Gamma_3(u_0)\cap\Gamma(u_2)$. The possible stabiliser subgroups $G_{u_0}$ give us two possibilities for the parameters $a_2, b_2$ as follows:

(a) if $G_{u_0}=PSL(2,9)$ or $P\Sigma L(2,9)$ then  $G_{u_0u_2}$ has two orbits in
$\Gamma(u_2)\setminus\{u_1,u_3\}$, each of length 4, and so $a_2=
|\Gamma_2(u_0)\cap\Gamma(u_2)|=4$ and $b_2=|\Gamma_3(u_0)\cap\Gamma(u_2)|=4$;

(b) in all other cases (namely $G_{u_0}=PGL(2,9), M_{10}$ or $P\Gamma L(2,9)$), $G_{u_0u_2}$ is
transitive on  $\Gamma(u_2)\setminus\{u_1,u_3\}$, and we have $a_2=|\Gamma_2(u_0)\cap\Gamma(u_2)|=0$ and $b_2=|\Gamma_3(u_0)\cap\Gamma(u_2)|=8$.
\qed

\medskip
We now show that case (2)(a) of Lemma~\ref{cover-g-1A} leads to no examples.

\begin{lemma}\label{cover-g-01}
There are no  $(G,3)$-geodesic-transitive graphs $\Gamma$ satisfying the conditions of
Lemma~$\ref{cover-g-1A}(2)(a)$.

\end{lemma}

\proof Let $\Gamma, \Sigma, G$ be as in Lemma~\ref{cover-g-1A} and suppose that part (2)(a) holds
so that $(a_2,b_2, c_2)= (4, 4,2)$ and $\Gamma$ has valency 10.


Let  $B,B'\in V(\Sigma)$ such that  $d_{\Sigma}(B,B')=2$, and let $u\in B$ and $u'\in B'$ such that $d_\Gamma(u,u')=2$. In our case (Lemma~\ref{cover-g-1A}(2)(a)), $G_u\cong G_B/N= PSL(2,9)$ or $P
\Sigma L(2,9)$, and $\Delta:=\Gamma(u')\cap \Gamma_2(u)$ is an orbit of $G_{uu'}$ of length 4.
Let $\Omega$ be the set of four blocks of $\Sigma_2(B)\cap\Sigma(B')$ containing a point of
$\Delta$, and let $\Omega'$ be the set consisting of the remaining four blocks of $\Sigma_2(B)\cap\Sigma(B')$. Then $G_{uu'}$ and $G_{BB'}$ leave $\Sigma_2(B)\cap\Sigma(B')=\Omega\cup\Omega'$ invariant, and both are transitive on $\Omega$. Moreover, the set $\Omega$ is the set of blocks adjacent to $B'$ in a $(G_B/N)$-invariant subgraph of  valency 4 of the induced subgraph
$[\Sigma_2(B)]$.

We consider the action of $G_{BB'}/N\leq \Aut(\Sigma)$ on $\Sigma_2(B)$. Now $\Aut(\Sigma)=PSL(3,4).\mathbb{Z}_2^2$ has a subgroup $Y$ such that $|Y:(G/N)|=2$ and
$Y_B=PGL(2,9)$ or $P\Gamma L(2,9)$ according as $G_u= PSL(2,9)$ or $P
\Sigma L(2,9)$, respectively. The group $Y_B$ is transitive on $\Sigma_2(B)$ of degree 45 and
the stabiliser $Y_{BB'}$ is a Sylow $2$-subgroup of $Y_B$. Similarly its index 2 subgroup $G_B/N$
is transitive  on $\Sigma_2(B)$ and its stabiliser $G_{BB'}/N$ is a Sylow $2$-subgroup of $G_B/N$. Thus each of $Y_B$ and $G_{B}/N$ has a unique transitive representation of degree 45 up to permutational isomorphism. A computation using MAGMA \cite{Magma-1997} shows that, in this representation, $Y_{BB'}$ has a unique orbit of size 8, and this orbit must be
$\Sigma_2(B)\cap\Sigma(B')$. A further computation shows that this $Y_{BB'}$-orbit
is the union of two orbits of $G_{BB'}/N$ of length 4.

Now $Y_B$ has a transitive action on the 45 flags (incident point-line pairs) of the
generalised quadrangle $GQ(2,2)$, and the uniqueness of this representation discussed above means that we may identify the vertices of $[\Sigma_2(B)]$ with the flags of $GQ(2,2)$.
Let us define two flags $(p,L)$ and $(q,M)$ (where $p,q$ are points and $L,M$ are lines of $GQ(2,2)$ incident with $p, q$, respectively) to be adjacent if and only if $p\ne q$, $L\ne M$, and either $p$ lies on $M$ or $q$ lies on $L$. Then $(p,L)$ is adjacent to exactly eight flags and $Y_B$ is transitive on ordered pairs of adjacent flags.
Thus this definition of adjacency defines an arc-transitive graph on $[\Sigma_2(B)]$ of valency 8,
and the uniqueness of $\Sigma_2(B)\cap\Sigma(B')$ as a $Y_{BB'}$-orbit of size 8 implies that the induced subgraph $[\Sigma_2(B)]$ has precisely this adjacency rule.  Further, taking $B'$ to be the flag $(p,L)$, one may check that
$G_{BB'}/N$ leaves invariant, and acts transitively on the sets $\{(q,M)\mid q\ne p, M\ne L, q\in L\}$ and
$\{(q,M)\mid q\ne p, M\ne L, p\in M\}$, each of size 4. It follows that $\Omega$ is equal to one of these sets. However the sets
$\{(q,M)\mid q\ne p, M\ne L, q\in L\}$ and
$\{(q,M)\mid q\ne p, M\ne L, p\in M\}$ are paired orbits of $G_{BB'}/N$, and therefore do not correspond to undirected graphs of valency 4. (Instead they correspond to $(G_{B}/N)$-invariant orientations of the edge set of $[\Sigma_2(B)]$.) This is a contradiction, since when $a_2=4$, the set $\Omega$ is the neighbourhood of $B'$ in a $(G_B/N)$-invariant valency 4 (undirected) subgraph of  $[\Sigma_2(B')]$.
Therefore $(a_2,c_2)\neq (4,2)$ for $\Gamma$ and there are no examples satisfying the conditions of Lemma~\ref{cover-g-1A}(2)(a).
\qed

\medskip
Finally we deal with case (2)(b) of Lemma~\ref{cover-g-1A}.

\begin{lemma}\label{cover-g-1}
Let $\Gamma$ be a $(G,3)$-geodesic-transitive graph, and suppose that the conditions of
Lemma~$\ref{cover-g-1A}(2)(b)$ hold. Then either
 \begin{itemize}
\item[(1)]  $\Gamma$ is the standard double cover of the Gewirtz graph, or

\item[(2)] $\Gamma$ has diameter at least $4$ and is not   $(G,4)$-distance-transitive.

\end{itemize}

\end{lemma}

\proof
First we show that the standard double cover $\overline{\Sigma}$ of the Gewirtz graph $\Sigma$ is an example.
Let $G=A\times\langle\tau\rangle$, where $A=\Aut(\Sigma)$ and $\tau$ is
the map defined before Lemma~\ref{rect-cover-1}. Then since $\Sigma$ is connected and non-bipartite with $c_2=2$,
it follows from  Lemma~\ref{rect-cover-1} that  $\overline{\Sigma}$ is  connected  with $c_2=2$. Also it follows from the definition that
$\overline{\Sigma}$ is a $G$-normal cover of $\overline{\Sigma}_{\overline N} \cong \Sigma$, where ${\overline N}= 1\times \langle\tau\rangle\unlhd G$.
Finally, it follows from the discussion at the end of the proof of Lemma~\ref{cover-g-1A} that $G$ is transitive on the $s$-geodesics
of  $\overline{\Sigma}$, for each $s\leq 3$. Thus all conditions hold and  $\overline{\Sigma}$ is an example.

If $\Gamma, G$ satisfy the hypotheses of the lemma, and if $\Gamma$ has diameter at least 4 and is not
$(G,4)$-distance-transitive, then (2) holds. So we assume that  this is not the case, and hence we assume that $\Gamma, G$ have the following properties:
 $\Gamma$ is connected,  $(G,3)$-geodesic-transitive, and $(G,s)$-distance-transitive for $s=\min\{4, \diam(\Gamma)\}$,
 $\Gamma$ has girth $4$ with $(a_2, b_2, c_2)=(0, 8, 2)$, $\Gamma_N\cong \Sigma$ for some non-trivial normal subgroup $N$ of $G$,
 and $G_{u_0}=PGL(2,9), M_{10}$ or $P\Gamma L(2,9)$, for $u_0\in V(\Gamma)$.

We claim that, if $\Gamma$ is not bipartite, then its standard double cover $\overline{\Gamma}$ and the
group $\overline{G}=G\times\langle\tau\rangle$ (with $\tau$ as defined before  Lemma~\ref{rect-cover-1}),
satisfy all of these conditions.
To prove this claim, assume that $\Gamma$ is not bipartite. Then $\overline{\Gamma}$ is  connected.
Since $\Gamma$ is connected of girth $4$ with $c_2=2$,
it follows from  Lemma~\ref{rect-cover-1} that  $\overline{\Gamma}$ has girth $4$  with $c_2=2$, and also from Lemma~\ref{rect-cover-1} it follows that
 $\overline{\Gamma}$ is $(\overline{G},3)$-geodesic-transitive (since $a_2=0$) and
 $\overline{\Gamma}$ is $(\overline{G},4)$-distance-transitive (since $\diam(\overline{\Gamma})>\diam(\Gamma)\geq3$).
 Finally $\overline{\Gamma}$ is a $\overline{G}$-normal cover of $\overline{\Gamma}_{\overline N} \cong \Gamma_N\cong \Sigma$,
 where ${\overline N}= N\times \langle\tau\rangle\unlhd \overline{G}$, and the vertex stabilisers in $\overline{G}$ and $G$ are isomorphic, implying, by Lemma~\ref{cover-g-1A} that  $(a_2, b_2)=(0, 8)$.
 This proves the claim.

 We now assume that $\Gamma$ is bipartite, if necessary replacing a non-bipartite graph with its standard double cover.
 We will prove that, under this assumption,  $\Gamma$ is $\overline{\Sigma}$ as in part (1).
 If in fact the original graph had been non-bipartite, then this would imply that  $\overline{\Gamma}=\overline{\Sigma}$ and hence that
 $\Gamma=\Sigma$, which is not the case since $\diam(\Sigma)=2$. Hence to  complete the proof of the lemma, it is sufficient to assume that
$\Gamma$ is bipartite and to prove that  $\Gamma=\overline{\Sigma}$.   So we assume that $\Gamma, G, N$ satisfy all the conditions
of the previous paragraph and in addition that $\Gamma$ is bipartite.

Note that  the Gewirtz graph $\Sigma$ has valency $10$, so also $\Gamma$ has valency 10 since it covers $\Sigma$.
As in the proof of  Lemma~\ref{cover-g-1A},
let  $(B_0,B_1,B_2)$ be a 2-geodesic of $\Sigma$ and let $B_3\in \Sigma(B_0)\cap \Sigma(B_2)$
such that $B_3\neq B_1$. Let $u_i\in B_i$, for each $i$, such that $(u_0,u_1,u_2,u_3)$ is a 3-arc of $\Gamma$; it was shown that this is a 4-cycle and $\Gamma(u_0)\cap \Gamma(u_2)=\{u_1, u_3\}$.
Let  $B_4\in \Sigma_2(B_0)\cap \Sigma(B_2)$ and let $u_4\in B_4$  so that $(u_0,u_1,u_2,u_4)$ is  a 3-arc; it was shown that this is a $3$-geodesic for
each of the $8$ choices for $B_4$ (since $a_2=8$).
Also in this case the stabiliser $G_{u_0}$ is $3$-transitive on  $\Gamma(u_0)$.

Let  $\Sigma(B_0)\cap \Sigma(B_4)=\{B_5,B_6\}$ and, for $i\in\{5,6\}$,
let $u_i\in B_i$ such that the sequence $(u_0,u_1,u_2,u_4,u_i)$ is a 4-arc.
Since $(u_0,u_1,u_2,u_4)$ is  a 3-geodesic, $d_\Gamma(u_0,u_i)\in\{2,3,4\}$, but  since
 $\Gamma$ is $(G,3)$-geodesic-transitive and $B_2, B_3\in\Sigma_2(B_0)$, all vertices of
 $ \Gamma_2(u_0)\cup\Gamma_3(u_0)$ lie in blocks of $\Sigma_2(B_0)$.
 Hence $u_i\in \Gamma_4(u_0)$, and since $\Gamma$ is $(G,4)$-distance-transitive,
 all vertices of  $ \Gamma_4(u_0)$ lie in blocks of $\Sigma(B_0)$. Now $\Gamma_4(u_0)\cap \Gamma(u_4)$ contains $\{u_5, u_6\}$ so $b_3\geq 2$.
Since $B_5\in \Sigma(B_0)$, $u_5$ is adjacent to a vertex $u_0'\in B_0$, and so
 $(u_0,u_1,u_2,u_4,u_5,u_0')$ is  a 5-arc. Thus $d_\Gamma(u_0,u_0')$ is at most 5, and is greater than $4$ since vertices of $\cup_{i=1}^4\Gamma_i(u_0)$ lie in blocks in
 $\Sigma(B_0)\cup \Sigma_2(B_0)$.
 Hence $d_\Gamma(u_0,u_0')=5$ and $\Gamma$ has diameter at least 5.
 Further, each vertex of $\Gamma_4(u_0)\cap\Gamma(u_4)$ lies in one of the two
blocks $B_5, B_6\in\Sigma(B_0)$ adjacent to $B_4$.
Thus  $\Gamma_4(u_0)\cap \Gamma(u_4)=\{u_5,u_6\}$ (since $\Gamma$ is a cover of $\Sigma$), and $b_3=2$.
Since $\Gamma$ is bipartite of valency 10, $a_3=0$ and  $c_3= |\Gamma_2(u_0)\cap\Gamma(u_4)|=   10-b_3=8$.
Thus we have $|\Gamma_2(u_0)|=\frac{10\times b_1}{c_2}=\frac{90}{2}=45$,
 $|\Gamma_3(u_0)|=\frac{45\times b_2}{c_3}=\frac{45\times 8}{8}=45$, and
 $|\Gamma_4(u_0)|=\frac{45\times b_3}{c_4}=\frac{90}{c_4}$, with $c_4\geq c_3=8$.
 Since vertices of $\Gamma_4(u_0)$ lie in blocks of $\Sigma(B_0)$,  $\Gamma_4(u_0)$ consists of a constant number of vertices from each of these blocks and so $|\Gamma_4(u_0)|$ is divisible by 10.
   It follows that $c_4=9$ and $|\Gamma_4(u_0)|=10$. Then $b_4=1$ (as $\Gamma$ is bipartite),
$\Gamma$ is $G$-distance-transitive (Remark \ref{3gt-dt-1}), and
 $|\Gamma_5(u_0)|=\frac{10\times b_4}{c_5}=\frac{10}{c_5}$, with $c_5\geq c_4=9$,
 so $c_5=10$ and $|\Gamma_5(u_0)|=1$.
Thus $\Gamma$ is an antipodal double cover of $\Sigma$ as well as bipartite, and has odd diameter.
Hence by Lemma~\ref{sdc}, $\Gamma$ is the standard double cover of the Gewirtz graph $\Sigma$.
 \qed

\bigskip
\bigskip

\bigskip
\subsection{Covers of the $M_{22}$ graph}

First we determine the value of $c_2$ for a $G$-normal cover of the $M_{22}$-graph.

\begin{lemma}\label{cover-m22-2}
Let $\Gamma$ be a  connected $(G,3)$-geodesic-transitive graph of girth $4$ or $5$.
Suppose that $\Gamma$ is a $G$-normal cover of the $M_{22}$-graph $\Sigma$, that is,  $G$ has a normal subgroup
$N$ such that $\Gamma_N\cong \Sigma$.
Let  $u\in V(\Gamma)$. Then
\begin{enumerate}
\item[(1)] $G/N=M_{22}$ or $M_{22}.2$, and  $G_u=\mathbb{Z}_2^4:S_6$ or $\mathbb{Z}_2^4:A_6$, respectively, and

\item[(2)] $c_2=4$.
\end{enumerate}

\end{lemma}
\proof (1) The $M_{22}$-graph $\Sigma$ is
strongly regular with parameters $(77, 16, 0, 4)$. Let $u\in B\in V(\Sigma)$.
By \cite{Brouwer-2018}, $A:=\Aut(\Sigma)\cong M_{22}.2$
and $A_B\cong \mathbb{Z}_2^4:S_6$. Moreover, it follows from \cite[p.39]{Atlas} that $M_{22}$
is the only proper subgroup of $M_{22}.2$ which is vertex transitive on $\Sigma$. Hence
$G/N= M_{22}$ or $M_{22}.2$, and $G_B/N= \mathbb{Z}_2^4:S_6$ or $\mathbb{Z}_2^4:A_6$, respectively.
Since $\Gamma$ is    $(G,3)$-geodesic-transitive  of girth $4$ or $5$, it follows from  Lemma \ref{cover-group-1} (2) that  $G_{B}=NG_u$, and
by Lemma \ref{2gt-covergt-1}, $N$ is semiregular on $V(\Gamma)$. Thus $G_u\cap N=1_G$ and
$G_u\cong G_u/(G_u\cap N) \cong NG_u/N= G_B/N$, proving part (1).

(2)  Since $\Gamma$ covers $\Sigma$ it follows that  $\Gamma$ has valency 16 and $c_2\leq 4$.
Suppose that $c_2<4$. Let  $(B_0,B_1,B_2)$ be a 2-geodesic of $\Sigma$. Let $u_i\in B_i$ ($0\leq i\leq 2$) such that
$(u_0,u_1,u_2)$ is a 2-arc of $\Gamma$. Since  $\Gamma$ has no triangles,
$(u_0,u_1,u_2)$ is a 2-geodesic.
Let $B_3\in \Sigma(B_0)\cap \Sigma(B_2)$ be one of the
$4-c_2$ blocks that contain no point of $\Gamma(u_0)\cap \Gamma(u_2)$,
and let $B_4$ be one of the 12 blocks of $\Sigma_2(B_0)\cap \Sigma(B_2)$.
Let $u_i\in B_i$ ($i=3, 4$) such that $(u_0,u_1,u_2,u_3)$
and $(u_0,u_1,u_2,u_4)$ are 3-arcs of $\Gamma$.
By Lemma~\ref{cover-group-1}(3), $u_3\in\Gamma_3(u_0)$ so that
$b_2:=|\Gamma(u_2)\cap \Gamma_3(u_0)|\geq 4-c_2$; and also all vertices
of $\Gamma_3(u_0)$ lie in blocks of $\Sigma(B_0)$ so that $u_4\not\in\Gamma_3(u_0)$.
 Since the block $B_4$ is not in $\Sigma(B_0)$ it follows that $u_4\in\Gamma_2(u_0)$,
 and hence that $a_2:=|\Gamma(u_2)\cap \Gamma_2(u_0)|\geq 12$. Since the
 valency $16=a_2+b_2+c_2$ we conclude that $a_2=12$ and $b_2=4-c_2$.

Since $B_3\in \Sigma(B_0)$, there exists a vertex  $u_0'\in B_0$ adjacent to $u_3$.
Further $u_0'\neq u_0$, as $u_3\not\in\Gamma(u_0)$.
 By  the $(G,3)$-geodesic-transitivity of  $\Gamma$, $G_{u_0}$ is transitive on $\Gamma_i(u_0)$ for $i=1,2,3$ and so vertices in these sets lie in blocks of $\Sigma(B_0)\cup\Sigma_2(B_0)$.
Hence $u_0'\in \Gamma_4(u_0)$, and $\Gamma$ has diameter at least 4.

Next we claim that $c_3:=|\Gamma(u_3)\cap\Gamma_2(u_0)|=15$. Consider
$\Delta:=\Gamma(u_3)\cap\Gamma_2(u_1)$. Since $d_\Gamma(u_1, u_3)=2$, $|\Delta|=a_2=12$. Since there are no $\Gamma$-edges between blocks of $\Sigma(B_0)$ and since all vertices of $\Gamma_3(u_0)$
lie in bocks of $\Sigma(B_0)$, it follows that $\Delta$ is disjoint from $\Gamma_3(u_0)$.
Then since $\Delta\subseteq \Gamma_2(u_1)$ it follows that $\Delta\subseteq\Gamma_2(u_0)$. Thus
$\Gamma(u_3)\cap\Gamma_2(u_0)$ contains $\Delta$ and so $c_3\geq 12$. To determine $c_3$ exactly we use some divisibility arguments. Since all vertices of $\Gamma_2(u_0)$ lie in blocks of $\Sigma_2(B_0)$, and since $G_{u_0}$ is transitive on both $\Gamma_2(u_0)$ and $\Sigma_2(B_0)$ (Lemma~\ref{cover-group-1}), it follows that $|\Sigma_2(B_0)|=60$ divides $|\Gamma_2(u_0)|=\frac{16\times 15}{c_2}=\frac{240}{c_2}$.
Since $c_2<4$, this implies that $c_2=1$ or 2. Again, by Lemma~\ref{cover-group-1}(3), $|\Gamma(u_0)|=16$ divides $|\Gamma_3(u_0)|= \frac{240}{c_2}\times\frac{4-c_2}{c_3}$, and this implies that $c_3$ divides $45$. It follows that $c_3=15$, as claimed.

Since $\Gamma$ has diameter at least 4, $b_3:=|\Gamma(u_3)\cap\Gamma_4(u_0)|=1$ and
$a_3:=|\Gamma(u_3)\cap\Gamma_3(u_0)|=0$,  and so by Remark \ref{3gt-dt-1}, $\Gamma$ is $G$-distance-transitive.
Since $a_2=12$, $\Gamma$ contains $5$-cycles and hence is  not bipartite,  and since $\Gamma$ is also not vertex-primitive, it follows from \cite[Theorem 2]{smith} that $\Gamma$ is antipodal.
We showed above that $c_2\leq 2$.
If $c_2=2$ then $|\Gamma_2(u_0)|=\frac{16\times 15}{c_2}=120$, $|\Gamma_3(u_0)|=\frac{240\times 2}{2\times c_3}=16$, and $|\Gamma_4(u_0)|=\frac{16\times 1}{c_4}$ with $c_4\geq c_3=15$, so we find that $|\Gamma_4(u_0)|=1$, $\Gamma$ has diameter 4,
and $\Gamma$ has intersection array $(16,15,2,1;1,2,15,16)$. However, by \cite[p.421]{BCN}, such a graph does not exist.
Hence $c_2=1$, so $b_2=3$ and we obtain  $|\Gamma_2(u_0)|=\frac{16\times 15}{c_2}=240$,
$|\Gamma_3(u_0)|=\frac{240\times 3}{c_3}=48$, and $|\Gamma_4(u_0)|=\frac{48\times 1}{c_4}$ with $c_4\geq c_3=15$. Thus $c_4=16$, $|\Gamma_4(u_0)|=3$. Hence $G_{u_0}$ has a transitive action on $\Gamma_4(u_0)$ of degree $3$, but this is impossible since by part (1), $G_{u_0}=
\mathbb{Z}_2^4:S_6$ or $\mathbb{Z}_2^4:A_6$.
\qed

\begin{lemma}\label{cover-m22-c3}
Let $\Gamma$ be a connected  $(G,3)$-geodesic-transitive  graph of girth $4$ or $5$.
Suppose that    $\Gamma$ is a $G$-normal cover of the $M_{22}$-graph. Then
either
 \begin{itemize}
\item[(1)]  $\Gamma$ is the standard double cover of the $M_{22}$-graph, or

\item[(2)] $\Gamma$ has diameter at least $4$ and is not   $(G,4)$-distance-transitive.

\end{itemize}
\end{lemma}

\proof  Let $\Sigma$ be the  $M_{22}$-graph, a strongly regular graph with parameters $(77, 16, 0, 4)$, and let $N\unlhd G$ such that $\Gamma_N\cong \Sigma$. Since   $\Gamma$ is a cover of  $\Sigma$, $\Gamma$ has valency 16, and by  Lemma \ref{cover-m22-2}, $G/N=M_{22}$ or $M_{22}.2$, and $c_2=4$. We identify $\Gamma_N$ with $\Sigma$.

Let $(u_0,u_1,u_2,u_3)$ be  a 3-geodesic of $\Gamma$ such that $u_i\in B_i\in V(\Sigma)$. Then   $(B_0,B_1,B_2,$ $B_3)$ is  a 3-arc of $\Sigma$.
Since $\Sigma$ also has $c_2=4$, it follows that   for each $B\in \Sigma(B_0)\cap \Sigma(B_2)$ such that $B\neq B_1$, the unique vertex of $\Gamma(u_2)\cap B$ is adjacent to  $u_0$. Hence
$B_3\notin \Sigma(B_0)$, and so  $B_3\in \Sigma_2(B_0)$.

Also  by   Lemma~\ref{cover-group-1}(4), we have   $G_{B_0B_2}=NG_{u_0B_2}
=N G_{u_0u_2}$ (since $c_2=4$).
It can be  easily checked using MAGMA \cite{Magma-1997} that $G_{B_0B_2}/N
\cong \mathbb{Z}_2^2.S_4$ or $\mathbb{Z}_2^2.S_4.\mathbb{Z}_2$ and acts transitively
on $\Sigma_2(B_0)\cap \Sigma(B_2)$, a set of 12 blocks including $B_3$. It follows that $G_{u_0u_2}$ acts transitively on $\Sigma_2(B_0)\cap \Sigma(B_2)$. Therefore $\Gamma_3(u_0)\cap \Gamma(u_2)$ contains a point from each of these 12 blocks, and hence $b_2:=|\Gamma_3(u_0)\cap \Gamma(u_2)|\geq 12$. Since the valency $16=a_2+b_2+c_2\geq a_2+12+4$, we have $a_2=0$ and $b_2=12$.

Now  $|\Sigma(B_0)\cap \Sigma(B_3)|=4$, say  $\Sigma(B_0)\cap \Sigma(B_3)=
\{E_1,E_2,E_3,E_4\}$, and, for each $i$, let
$ \Gamma(u_3)\cap E_i=\{e_i\}$. Then $d_\Gamma(u_0,e_i)\leq 4$, and $d_\Gamma(u_0,e_i)\ne 1$ since $u_3\in\Gamma_3(u_0)$.  Since $G_{u_0}$ is transitive on $\Gamma_i(u_0)$ for $i=2, 3$ and $u_2, u_3$ lie in blocks in $\Sigma_2(B_0)$, it follows that $d_\Gamma(u_0,e_i)\ne 2, 3$ and hence $d_\Gamma(u_0,e_i)= 4$, so $\diam(\Gamma)\geq4$ and $b_3=|\Gamma_4(u_0)\cap\Gamma(u_3)|\geq 4$.

If $\Gamma$ is not $(G,4)$-distance-transitive then part (2) holds, so suppose now that $\Gamma$ is   $(G,4)$-distance-transitive.
Then all vertices in $\Gamma_4(u_0)$ lie in blocks of $\Sigma(B_0)$, and hence
vertices in $B_0\setminus\{u_0\}$ all have distance at least $5$ from $u_0$. Since $\Gamma$ covers $\Sigma$ there is a unique vertex $u_0'\in B_0\setminus\{u_0\}$ adjacent to $e_1$ and this must satisfy $d_\Gamma(u_0,u_0')= 5$. 
Further, since each vertex $u\in\Gamma_4(u_0)\cap \Gamma(u_3)$ lies in one of the 4 blocks of $\Sigma(B_0)\cap \Sigma(B_3)$, it follows that $u$ is one of the $e_i$, and so $b_3=|\Gamma_4(u_0)\cap \Gamma(u_3)|=4$.

Suppose that $\Gamma$ is not  bipartite, and let $\overline{\Gamma}$ denote the  standard double cover  of $\Gamma$ and $\overline{G}=G\times \mathbb{Z}_2$, as in Lemma \ref{rect-cover-1}. Then, by  Lemma \ref{rect-cover-1},
$\overline{\Gamma}$ is $(\overline{G},3)$-geodesic-transitive and $(\overline{G},4)$-distance-transitive and has $c_2=4$. Moreover
$\overline{\Gamma}_{\overline N}\cong \Sigma$ where ${\overline N}=N\times \mathbb{Z}_2$.  If Lemma~\ref{cover-m22-c3} holds for bipartite graphs, then $\overline{\Gamma}$ is the standard double cover of $\Sigma$, and hence $\Gamma=\Sigma$, which is a contradiction.
Thus it is sufficient to prove the lemma for bipartite graphs, and we therefore assume from now on that $\Gamma$ is bipartite.

Then, since $\Gamma$ is bipartite, $a_3=|\Gamma_3(u_0)\cap\Gamma(u_3)|=0$,
so $c_3=|\Gamma_2(u_0)\cap\Gamma(u_3)|=16-a_3-b_3=12$, $|\Gamma_3(u_0)|=\frac{60\times b_2}{c_3}=60$, and $|\Gamma_4(u_0)|=\frac{60\times b_3}{c_4}=\frac{240}{c_4}$ with $c_4\geq c_3=12$.
Since $\Gamma_4(u_0)$ is contained in the union of the blocks in $\Sigma(B_0)$, and since $G_{u_0}$ is transitive on $\Sigma(B_0)$ and fixes
$\Gamma_4(u_0)$ setwise, it follows that $\Gamma_4(u_0)$ contains a constant number of vertices from each block of $\Sigma(B_0)$. Thus $|\Gamma_4(u_0)|$ is divisible by $16$, and it follows that $c_4=15$ and $|\Gamma_4(u_0)|=16$. Again, since $\Gamma$ is bipartite,  $a_4=|\Gamma_4(u_0)\cap\Gamma(e_1)|=0$
so $b_4=|\Gamma_5(u_0)\cap\Gamma(e_1)|=16-c_4=1$. It follows that $G_{u_0}$ is transitive on
$\Gamma_5(u_0)$ and that $|\Gamma_5(u_0)|=\frac{16\times b_4}{c_5}=\frac{16}{c_5}$ with $c_5
= |\Gamma_4(u_0)\cap\Gamma(u_0')|\geq c_4=15$. Thus $c_5=16$, $|\Gamma_5(u_0)|=1$, and $\Gamma$
is a $G$-distance-transitive, bipartite and antipodal graph, with antipodal blocks of size 2 and with odd diameter $5$. It follows from Lemma~\ref{sdc} that   $\Gamma$ is the standard double cover of $\Sigma$, as in part (1).
 \qed

\subsection{Covers of  complete bipartite graphs}


\begin{lemma}\label{cover-bipart-1}
Let $\Gamma$ be a connected  $(G,3)$-geodesic-transitive graph of girth $4$ or $5$.
Suppose that  $\Gamma$ is a $G$-normal $m$-fold cover of $\K_{r,r}$,  where $r\geq 3$ and $m\geq 2$. Then $\Gamma$ is bipartite with  girth $4$ and diameter at least $4$, and
one of the following holds:

 \begin{itemize}
\item[(1)]  $\Gamma$ is  a  Hadamard graph of order $r$ (so   $m=2$); or

\item[(2)] $\Gamma$ is a $G$-distance-transitive antipodal graph $m\K_{r,r}$ which is the incidence graph of a resolvable divisible design
$RGD(r,c_2,m)$, where  $r=mc_2$, such that any two blocks from different parallel classes contain exactly $c_2$ common points; $\Gamma$ has   intersection array $(r,r-1,r-c_2,1;1,c_2,r-1,r)$, and $c_2\geq 2$, $m\geq3$; or

\item[(3)]  $\Gamma$ is  not $(G,4)$-distance-transitive.

\end{itemize}

\end{lemma}

There are examples of distance-transitive antipodal graphs $m\K_{r,r}$ arising from resolvable divisible designs $RGD(r,c_2,m)$ in case (2). Those with valency  $r\leq 13$ are classified and are listed in  \cite[p.223--225]{BCN}. Moreover, by Lemma~\ref{res-divisible-1}, each $RGD(r, c_2, m)$ with $r=c_2m$, such that any two blocks from different parallel classes contain exactly $c_2$ common points, has an incidence graph which is distance-regular.

\proof  Let $\Sigma=\K_{r,r}$ where $r\geq 3$, and let $N\unlhd G$ be such that $\Gamma_{N}\cong \Sigma$.
Since  $\Sigma$ is bipartite it follows that  $\Gamma$ is also bipartite, and hence contains no cycles of odd length, Then since
$\Gamma$ has girth 4 or 5, the girth of  $\Gamma$ must be  4.
Set $V(\Sigma)=\mathcal{B}_0\cup \mathcal{B}_1$ where
$\mathcal{B}_0=\{B_{01},B_{02},\ldots,B_{0,r}\}$ and
$\mathcal{B}_1=\{B_{11},B_{12},\ldots,B_{1,r}\}$, and $\mathcal{B}_0$, $\mathcal{B}_1$ are  the two bipartite halves of $\Sigma$.

Let $(u_{01},u_{11},u_{02},u_{12})$ be a 4-cycle of $\Gamma$ such that each $u_{ij}\in B_{ij}$. Then $u_{02}\in  \Gamma_2(u_{01})$ and  $(B_{01},B_{11},B_{02},B_{12})$ is a 4-cycle of $\Sigma$.
Since  $\Gamma$ is $(G,3)$-geodesic-transitive with  diameter at least 3,   $u_{02}$ is  adjacent to some vertex  $u_{13}\in  \Gamma_3(u_{01})$. Let $u_{13}\in B_{13}$. 
Since $\Sigma=\K_{r,r}$, it follows that  $B_{13}\in \Sigma(B_{01})$.
Let $u_{01}'\in B_{01}$ be adjacent to  $u_{13}$. Then $d_\Gamma(u_{01}, u_{01}')\leq 4$. As $\Gamma$ is bipartite
 $d_\Gamma(u_{01}, u_{01}')$ is even, and as $\Gamma$ is a cover of $\Sigma$,  $d_\Gamma(u_{01}, u_{01}')\ne 2$. Hence
$u_{01}'\in  \Gamma_4(u_{01})$, and in particular   $\diam(\Gamma)\geq 4$.

Since $\Gamma$ has girth 4,  it follows that  $c_2\geq 2$, and we have $|\Gamma_2(u_{01})|=\frac{r\times (r-1)}{c_2}$.
By Lemma~\ref{2gt-covergt-1}, it follows that $G_{B_{01}}$ is transitive on $\Sigma_2(B_{01})=  \mathcal{B}_0\setminus\{B_{01}\}$,
and by Lemma~\ref{cover-group-1}, $G_{B_{01}}=NG_{u_{01}}$, and hence also $G_{u_{01}}$ is transitive on
 $\Sigma_2(B_{01})$. It follows that $\Gamma_2(u_{01})$ contains equally many vertices from each block of $\Sigma_2(B_{01})$ and hence
 that $r-1$ divides  $|\Gamma_2(u_{01})|$. Therefore $c_2$ divides $r$. Since  $\diam(\Gamma)\geq 4$, we must have $c_2<r$, and hence
 $2\leq c_2\leq r/2$.

Suppose first that  $c_2=\frac{r}{2}$. Then by  \cite[Theorem 1.9.3]{BCN} (since $\Gamma$ has diameter at least $4$ and valency $r\geq3$),
$\Gamma$ is a Hadamard graph of order $r$ with intersection array
$(r, r-1,r/2,1;1,r/2,r-1,r)$. Thus $|V(\Gamma)|=4r$ and so $m=|B_{01}|=2$ and (1) holds.

Now suppose that  $2\leq c_2< \frac{r}{2}$, so $c_2\leq r/3$ since $c_2$ divides $r$.
If  $\Gamma$ is not $(G,4)$-distance-transitive then (3) holds, so
assume now that $\Gamma$ is $(G,4)$-distance-transitive.
Then since  $u_{01}'\in \Gamma_4(u_{01})\cap B_{01}$, it follows that  $\Gamma_4(u_{01})\subset B_{01}$.
Let $u\in \Gamma(u_{01}')\setminus\{u_{13}\}$. Then $u\in \Gamma_3(u_{01})\cup \Gamma_5(u_{01})$.
Further, $(u_{13},u_{01}',u)$ is a 2-geodesic, so $|\Gamma(u_{13})\cap \Gamma(u)|=c_2\geq 2$, say
$\{u_{01}',u_{03}\}\subseteq \Gamma(u_{13})\cap \Gamma(u)$.
Since $u_{13}\in \Gamma_3(u_{01})$, it follows that $u_{03}\in \Gamma_2(u_{01})\cup \Gamma_4(u_{01})$. Moreover, since
$\Gamma_4(u_{01})\subset B_{01}$, we have $u_{03}\notin  \Gamma_4(u_{01})$, and so $u_{03}\in \Gamma_2(u_{01})$. Hence
$u\in \Gamma_3(u_{01})$, and it follows that $\Gamma(u_{01}')\subseteq \Gamma_3(u_{01})$ so $c_4=r$ and $\diam(\Gamma)=4$.
In particular $\Gamma$ is $G$-distance-transitive.  Since $\Gamma$ is a cover of $\Sigma$, each vertex of $B_{01}\setminus\{u_{01}\}$
has even distance from $u_{01}$ greater than 2, and therefore $B_{01}\setminus\{u_{01}\}=\Gamma_4(u_{01})$ so $\Gamma$ is an antipodal graph.
Also $\Gamma_2(u_{01})=\cup_{i=2}^rB_{0i}$ and $\Gamma_3(u_{01})=(\cup_{i=1}^rB_{1i})\setminus \Gamma(u_{01})$, so
$|\Gamma_i(u_{01})|$ is $(r-1)m, r(m-1), m-1$ for $i=2,3,4$, respectively. In particular, since
$|\Gamma_2(u_{01})|= \frac{r}{c_2}(r-1)$, it follows that $r=c_2m$ so $m\geq3$ since $c_2\leq r/3$. Also, counting the edges between $\Gamma_3(u_{01})$ and
$\Gamma_4(u_{01})$ we have $r(m-1)b_3=|\Gamma_3(u_{01})|b_3=|\Gamma_4(u_{01})|c_4=(m-1)r$ and hence $b_3=1$.
Therefore $c_3=r-b_3=r-1$ (since $a_3=0$),
and $\Gamma$ has  intersection array
\begin{center}
$(r,r-1,r-c_2,1;1,c_2,r-1,r)$.
\end{center}
Moreover, by \cite[Item 10 on p.316]{FII-1986}, a distance-regular graph  with intersection array of this type
is an $m$-fold antipodal cover $m\K_{r,r}$ where $r=m.c_2$.
Finally, we show that $\Gamma$ is the incidence graph of a resolvable divisible design with the properties specified in (2), and that each such design
determines a distance-regular  $m$-fold antipodal cover of $\K_{r,r}$.

Let $V(\Gamma)=V_1\cup V_2$ where $V_1,V_2$ are  the bipartite halves  of $V(\Gamma)$. We construct the design
$\mathcal{D}=(V_1,\mathcal{P},  \mathcal{B})$ with $V_1$ being the set of points, $\mathcal{P}$ being the set of antipodal blocks of $\Gamma$ contained in $V_1$ (that is, $\mathcal{P} = \mathcal{B}_0$), and $\mathcal{B}$
being the set $\{ \Gamma(u)\mid u\in V_2\}$ of neighbour-sets of vertices in $V_2$.
Then  $\mathcal{P}$ consists of $r$  antipodal blocks, each  of size $m$, and is a partition of $V_1$.
Since $\Gamma$ is a bipartite graph of valency $r$ that covers $\Sigma$, it follows that each block in $\mathcal{B}$  is adjacent to exactly one vertex of every antipodal block in $\mathcal{P}$. Also, each pair of points from different antipodal blocks of $\mathcal{P}$ are at distance $2$ in $\Gamma$ and therefore are contained in exactly $c_2$ blocks of $\mathcal{B}$.
Thus  $\mathcal{D}$ is a  divisible design $GD(r,c_2,m,rm)$.
Moreover, the set $\mathcal{B}$ can be partitioned into $r$ parts, namely $\mathcal{B}^{(i)}:=\{ \Gamma(u)\mid u\in B_{1i}\}$ for $1\leq i\leq r$.
Each of these parts has size $m$ and distinct $\Gamma(u), \Gamma(u')$ in the same part must be disjoint since $\Gamma$ is a cover of $\Sigma$.
Thus $\cup_{u\in B_{1i}}\Gamma(u)$ has size $rm$ and hence is equal to $V_1$, that is, each  $\mathcal{B}^{(i)}$ is a partition of $V_1$ -- a parallel class of blocks of $\mathcal{D}$. Thus $\mathcal{D}$  is a resolvable divisible design $RGD(r,c_2,m)$. If $\Gamma(u)\in \mathcal{B}^{(i)}$ and
$\Gamma(v)\in \mathcal{B}^{(j)}$ with $i\ne j$, then $u, v$ are at distance $2$ in $\Gamma$ and hence $|\Gamma(u)\cap\Gamma(v)|=c_2$. Thus $\mathcal{D}$ has all the properties specified in (2).

 Conversely, by Lemma \ref{res-divisible-1}, if $\mathcal{D}$ is a resolvable divisible
design $RGD(r,c_2,m)$, where $r=c_2.m$, such that any two blocks from different parallel classes contain exactly $c_2$ common points, then
its incidence graph $\Inc(\mathcal{D})$ is an  antipodal diameter 4 graph that is an $m$-fold cover $m\K_{r,r}$.
\qed

\medskip


\subsection{Covers of the Hoffman-Singleton graph}

\medskip
Recall that the automorphism  group of the Hoffman-Singleton graph is $PSU(3,5).\mathbb{Z}_2$. Its action on vertices has rank $3$, with point stabiliser $S_7$, see  \cite[p.304]{FIKW-1}.

\begin{lemma}\label{ho-s-2}
Let $\Gamma$ be a connected  $(G,3)$-geodesic-transitive  graph of girth $4$ or $5$.
Then   $\Gamma$  is not  a $G$-normal cover of the Hoffman-Singleton graph.

\end{lemma}
\proof Suppose that   $\Gamma$  is a $G$-normal cover of the Hoffman-Singleton graph $\Sigma$, and let $N\unlhd G$
such that  $\Gamma_{N}= \Sigma$.
Since $\Sigma$ is a strongly regular graph  with parameters $(50,7,0,1)$,
it follows that   $\Gamma$ has valency 7 and $c_2=1$. In particular, $\Gamma$   has  girth  5.

Let $(u_0,u_1,u_2,u_3)$ be  a 3-geodesic of $\Gamma$ such that $u_i\in B_i\in V(\Sigma)$. Then   $(B_0,B_1,B_2,$ $B_3)$ is  a 3-arc of $\Sigma$ and $|\Sigma_2(B_0)|=42$.
Since $\Sigma$  has $c_2=1$,   $B_3$ lies in $ \Sigma_2(B_0)$.
By Lemma~\ref{2gt-covergt-1},  $\Sigma$ is    $(G/N,2)$-geodesic-transitive,
so  $G_{B_0}/N$  is 2-transitive on $\Sigma(B_0)$. By \cite[p34]{Atlas}, the only proper subgroup of
$\Aut(\Sigma)$ that is $2$-geodesic-transitive is $PSU(3,5)$, and so
$G/N=PSU(3,5).\mathbb{Z}_2$ or $PSU(3,5)$ and
$G_{B_0}/N=S_7$ or $A_7$, respectively, acting faithfully on $\Sigma(B_0)$.
For both cases, we check using MAGMA \cite{Magma-1997} that  $G_{B_0B_2}/N$ is transitive on
$\Sigma_2(B_0)\cap \Sigma(B_2)$, a set of size $6$.
By   Lemma \ref{cover-group-1}(4),   $G_{B_0}=NG_{u_0}$ and $G_{B_0B_2}=NG_{u_0u_2}=NG_{u_0B_2}$.
Therefore,    $G_{u_0u_2}$ is  transitive on $\Sigma_2(B_0)\cap \Sigma(B_2)$.
Let $C \in  \Sigma_2(B_0)\cap \Sigma(B_2)$ with $C\ne B_3$, and let $ \Gamma(u_2)\cap C=\{c\}$.
Then some element $g\in G_{u_0u_2}$  maps  $B_3$ to  $C$, and
so $u_3^g=c$. Thus $c\in \Gamma_3(u_0)\cap \Gamma(u_2)$.
Since $|\Sigma_2(B_0)\cap \Sigma(B_2)|=6$, it follows that $|\Gamma_3(u_0)\cap \Gamma(u_2)|=6$, that is, $\Gamma$ has $b_2=6$, and so $a_2=0$,
contradicting the fact that  $\Gamma$   has  girth  5.  \qed

\medskip
We are ready to prove the third  theorem.

\medskip

\noindent {\bf Proof of Theorem \ref{girth45-redtheo-1}.}
Let $\Gamma$ be a connected  $(G,3)$-geodesic-transitive  graph of girth  4 or 5, and let
$N\unlhd G$ with at least $3$ orbits on $V(\Gamma)$, and such that  $\Gamma_N$ has  diameter at most $2$.
It follows from  Theorem \ref{3gt-redtheo-1} that
$\Gamma$ is a cover of $\Gamma_N$,    and    $\Gamma_N$ is either a complete graph or
a strongly regular graph of girth $4$ or $5$.
If $\Gamma_N$ is a complete graph, then by  Lemma \ref{quot-comp-2},
$\Gamma$ is isomorphic to  $\K_{k+1,k+1}-(k+1)\K_2$, for some $k\geq 3$, or  to  $[HoS]_2$, as in Table \ref{table-cover-diam2}.
Suppose that $\Gamma_N$ has  diameter $2$ and girth $4$ or $5$. Then by Theorem \ref{diam2-theo-1},
$\Gamma_N$  is one of the following graphs: $C_5$, $\K_{r,r}$ with $r\geq 2$, the Higman-Sims graph HiS,  the Gewirtz graph, the $M_{22}$-graph, the folded $5$-cube $\Box_5$,
the Petersen graph,  and the Hoffman-Singleton graph. Further, by Lemmas  \ref{cover-petersen}(1) and \ref{ho-s-2},
$\Gamma_N$  is not $C_5$ or the Hoffman-Singleton graph. All the other possibilities yield examples, and
by Lemmas  \ref{cover-petersen} (2)-(3), \ref{cover-hs-1}, \ref{cover-g-1}, \ref{cover-m22-c3} and  \ref{cover-bipart-1},  either $\Gamma$
is one of the graphs in Table~\ref{table-cover-diam2}, or $\Gamma$ is not $(G,4)$-distance-transitive, with $\Gamma_N$ as in Table~\ref{table-cover-diam}.
\qed


\begin{thebibliography}{hhhh}

\bibitem{ACMM-1996}
B. Alspach, M. Conder, D. Maru$\check{{\rm s}}$i$\check{{\rm c} }$ and M. Y. Xu, A
classification of 2-arc-transitive circulants, {\it J. Algebraic
Combin.} {\bf 5} (1996), 83--86.

\bibitem{Aschb-1971}
M. Aschbacher, The nonexistence of rank three permutation groups of degree 3250 and subgree 57, {\it J. Algebra}, {\bf 19} (1971), 538--540.



\bibitem{B-1972}
E. Bannai, Maximal subgroups of low rank of finite symmetric and alternating groups,
{\it J. Fac. Sci. Univ. Tokyo} {\bf 18} (1972), 475--486.

\bibitem{Biggs-1}
N. L. Biggs, {\it Algebraic Graph Theory}, Cambridge University Press, New
York, (1974).

\bibitem{Magma-1997}
W. Bosma, C. Cannon, and C. Playoust, The MAGMA algebra system I: The user
language, {\it J. Symbolic Comput.} {\bf 24} (1997), 235--265.


\bibitem{Brouwer-2018}
A. E. Brouwer, $M_{22}$ Graph, Technische Universiteit Eindhoven, \url{http://www.win. tue.nl/~aeb/graphs/M22.html}. Accessed 28 March, 2019.

\bibitem{BCN}
A. E. Brouwer, A. M. Cohen and A. Neumaier, {\it Distance-Regular Graphs},
Springer Verlag, Berlin, Heidelberg, New York, (1989).

\bibitem{BM-1994}
F. Buekenhout and H. Van Maldeghem, A characterization of some rank
2 incidence geometries by their automorphism group, {\it  Mitt. Math. Sem.
Giessen} {\bf 218} (1994), 1--70.





\bibitem{Cameron-1}
P. J. Cameron,  {\it Permutation Groups}, volume 45 of London Mathematical
Society Student Texts, Cambridge University Press, Cambridge,
(1999).

\bibitem{CP-1983}
P. J. Cameron and C. E. Praeger,  On 2-arc transitive graphs of girth 4, {\it J. Combin. Theory Ser. B} {\bf 35} (1983), 1--11.

\bibitem{Atlas}
J. H. Conway, R. T. Curtis, S. P. Norton, R. A. Parker and R. A.
Wilson, Atlas of Finite Groups,  {\it Clarendon Press} Oxford,
(1985).

\bibitem{DGLP-locdt-2012}
A. Devillers, M. Giudici, C. H. Li and C. E. Praeger, Locally
$s$-distance transitive Graphs, {\it J. Graph Theory } {\bf (2)69}
(2012), 176--197.


\bibitem{DJLP-clique}
A. Devillers, W. Jin, C. H. Li and C. E. Praeger, Local $2$-geodesic
transitivity and clique graphs, {\it J. Combin. Theory Ser. A} {\bf
120} (2013), 500--508.

\bibitem{DJLP-2}
A. Devillers, W. Jin, C. H. Li and C. E. Praeger, Line graphs and
geodesic transitivity, {\it Ars Math. Contemp.} {\bf 6} (2013),
13--20.


\bibitem{DJLP-cayleyred1}
A. Devillers, W. Jin, C. H. Li and C. E. Praeger, On normal
2-geodesic transitive Cayley graphs, {\it J. Algebr. Comb.} {\bf 39} (2014),
903--918.

\bibitem{DJLP-prime}
A. Devillers, W. Jin, C. H. Li and C. E. Praeger, Finite 2-geodesic
transitive graphs of prime valency, {\it J. Graph Theory} {\bf
80} (2015), 18--27.


\bibitem{DM-1}
J. D. Dixon and B. Mortimer,  {\it Permutation Groups}, Springer, New
York, (1996).




\bibitem{FII-1986}
I. A. Faradzev, I. A. Ivanov and  A. A. Ivanov, Distance-transitive graphs of valency 5, 6 and 7,  {\it European J. Combin.}
{\bf 7} (1986),  303--319.



\bibitem{FIKW-1}
I. A. Faradzev, A. A. Ivanov, M. H. Klin and A. J. Woldar, {\it Investigations in Algebric Theory of Combinatorial Objects}, Springer-Science+Business Media, B.V., (1994).

\bibitem{Frucht-1}
R. Frucht,   Die gruppe des Petersen'schen Graphen und der Kantensysteme der regularen Polyeder, {\it Comment. Math. Helv.}, {\bf 9} (1937), 217--223.




\bibitem{GR}
C. D. Godsil and G. Royle, {\it Algebraic Graph Theory}, Springer, New York,
Berlin, Heidelberg,  (2001).


\bibitem{GLP}
C. D. Godsil, R. A. Liebler and C. E. Praeger,  Antipodal distance
transitive covers of complete graphs, {\it European J. Combin.} {\bf
19} (1998), 455--478.



\bibitem{IP-1}
A. A. Ivanov and C. E. Praeger,  On finite affine 2-arc transitive
graphs, {\it European J. Combin.} {\bf 14} (1993), 421--444.

\bibitem{WJ-2015au}
W. Jin,  Finite 3-geodesic transitive but not 3-arc transitive graphs, {\it Bull. Aust. Math. Soc.} {\bf 91} (2015), 183--190.

\bibitem{JW-2018}
W. Jin, The pentavalent three-geodesic-transitive graphs, {\it Discrete Math.} {\bf 341} (2018), 1344--1349.


\bibitem{KL-rank3-1982}
W. M. Kantor and R. A. Liebler, The rank 3 permutation
representations of the finite classical groups, {\it Trans. Amer.
Math. Soc.} {\bf 271} (1982), 1-71.



\bibitem{Li-abeliancay-2008}
C. H. Li and J. M. Pan, Finite 2-arc-transitive abelian Cayley
graphs, {\it European J. Combin. } {\bf 29} (2008), 148--158.


\bibitem{LPS-1}
M. W. Liebeck, C. E. Praeger and J. Saxl,  On the O'Nan-Scott
theorem for finite primitive permutation groups,  {\it J. Austral.
Math. Soc. Ser. A} {\bf 44} (1988), 389--396.

\bibitem{LS-1986}
M. W. Liebeck and J. Saxl, The finite primitive permutation groups
of rank three, {\it Bull. London Math. Soc.} {\bf 18(2)} (1986),
165--172.




\bibitem{Praeger-1988-2trans}
C. E. Praeger, Primitive permutation groups with a doubly transitive subconstituent, {\it
J. Austral Math. Soc. (Series A)} {\bf 45} (1988) 66--77.







\bibitem{PraegerWang-1996}
C. E. Praeger and J. Wang, On primitive representations of finite alternating and symmetric groups with a 2-transitive subconstituent, {\it
J. Algebra} {\bf 180} (1996) 808--833.


\bibitem{smith}
D. H. Smith, Primitive and imprimitive graphs, {\it Quart. J. Math.
Oxford} {\bf (2)22} (1971), 551--557.


\bibitem{Tutte-1}
W. T. Tutte,  A family of  cubical  graphs, {\it Proc. Cambridge
Philos. Soc.} {\bf 43} (1947), 459--474.


\bibitem{Tutte-2}
W. T. Tutte, On the symmetry of cubic graphs, {\it Canad. J. Math.}
{\bf 11} (1959), 621--624.




\bibitem{weiss}
R. Weiss,  The non-existence of 8-transitive graphs, {\it
Combinatorica} {\bf 1} (1981), 309--311.

\bibitem{Wielandt-book}
H. Wielandt, {\it Finite Permutation Groups}, New York: Academic Press
(1964).


\end{thebibliography}
\end{document}